\documentclass[reqno,11pt]{amsart}
\usepackage{amssymb, amsmath,latexsym,amsfonts,amsbsy, amsthm}
\usepackage{color,esint,pgfplots}
\usepackage{cases}
\usepackage{mathtools}
\usepackage{cancel}

\usepackage{hyperref}
\hypersetup{hidelinks}

\usepackage{enumitem}
\usepackage{setspace}
\usepackage{mathrsfs}
\usepackage{stmaryrd}
\usepackage{graphicx}
\usepackage{comment}
\allowdisplaybreaks[4]
\usepackage[margin=1.2in]{geometry}

\newcommand{\picdis}[1]{}

 \newcommand{\eps}{\epsilon}

\newcommand{\vep}{\varepsilon}

\renewcommand{\[}{\begin{equation}}
\renewcommand{\]}{\end{equation}}

\renewcommand{\div}{\operatorname{div}}

\newcommand{\eqdef}{\overset{\mbox{\tiny{def}}}{=}}

\newtheorem{theorem}{Theorem}[section]

\newtheorem{lemma}[theorem]{Lemma}
\newtheorem{proposition}[theorem]{Proposition}

\theoremstyle{remark}
\newtheorem{remark}{Remark}[section]

\begin{document}
\title[Stability of viscoelastic boundary layers]{Effects of elasticity on the stability of viscoelastic boundary layers}

\author{Shengxin Li}
\address{School of Mathematics and Statistics, Ningbo University, Ningbo, China.}
\email{lishengxin@nbu.edu.cn}
\author{Dehua Wang}
\address{Department of Mathematics, University of Pittsburgh, Pittsburgh, USA.}
\email{dhwang@pitt.edu}
\author{Feng Xie}
\address{School of Mathematical and Sciences, Shanghai Jiao Tong University, Shanghai, China.}
\email{tzxief@sjtu.edu.cn}
\author{Zhu Zhang}
\address{Department of Applied Mathematics, The Hong Kong Polytechnic University, Hong Kong, China.}
\email{zhuama.zhang@polyu.edu.hk}

\begin{abstract}
 
The elastic effect on the stability of boundary layers is studied for incompressible elastic fluids in the half-plane. It is known from [33] that if the deformation tensor matrix is nonsingular near the boundary, elasticity can prevent the formation of a boundary layer. In this paper, we prove stability properties when the deformation tensor degenerates. Our results are twofold: (1) when the first row of the deformation tensor has a non-degenerate projection onto a fixed direction, the initial-boundary value problem for the two-dimensional viscoelastic Prandtl boundary layer equations in the half-plane admits local well-posedness in Sobolev spaces, without any monotonicity condition on the velocity profile; (2) when the first row of the deformation tensor vanishes simultaneously at one point, the linearized Prandtl boundary layer system is ill-posed. The first result reveals the stabilizing effect of elasticity against the classical Prandtl instability. The second result shows that elastic effects cannot suppress the Prandtl-type instability associated with a non-degenerate critical point of the velocity profile.

\end{abstract}

\date{\today}

\maketitle

%
%
%
\numberwithin{equation}{section}

\setcounter{secnumdepth}{3}

\indent


\section{Introduction}
In this paper, we consider the following initial-boundary value problem for the two-dimensional viscoelastic Prandtl boundary layer equations in the half-plane $\Omega=\{(x, y)\mid  x\in\mathbb{T}, y\in \mathbb{R}_+\}$:
\begin{align}\label{1.1}
\begin{cases}
\partial_tu+u\partial_xu+v\partial_yu-\partial_y^2u-f_1\partial_xf_1-f_3\partial_yf_1-f_2\partial_xf_2-f_4\partial_yf_2=-P(t,x),\\
\partial_tf_1+u\partial_xf_1+v\partial_yf_1-f_1\partial_xu-f_3\partial_yu=0,\\
\partial_tf_2+u\partial_xf_2+v\partial_yf_2-f_2\partial_xu-f_4\partial_yu=0,\\
\partial_xu+\partial_yv=0,\quad \partial_xf_1+\partial_yf_3=0,\quad \partial_xf_2+\partial_yf_4=0,\\
(u, v, f_3, f_4)|_{y=0}={\bf{0}},\quad \lim\limits_{y\to\infty}(u, f_1, f_2)(t, x, y)=(U, F_1, F_2)(t, x),\\
(u, f_1, f_2)|_{t=0}=(u_0, f_{1, 0}, f_{2, 0})(x, y).
\end{cases}
\end{align}
Here the vector $(u, v)$ and the matrix $\left( \begin{array}{cc} f_1& f_2\\ f_3&  f_4 \end{array}\right)$ represent the velocity field and the deformation tensor, respectively.  The pressure function $P(t,x)$ and the far fields functions  $(U, F_1, F_2)(t,x)$ are given, and satisfy the following Bernoulli law:
\begin{align}\label{1.2}
\begin{cases}
\partial_tU+U\partial_xU+\partial_xP-F_1\partial_xF_1-F_2\partial_xF_2=0,\\
\partial_tF_1+U\partial_xF_1-F_1\partial_xU=0,\\
\partial_tF_2+U\partial_xF_2-F_2\partial_xU=0.
\end{cases}
\end{align}
The system \eqref{1.1}-\eqref{1.2} is an extension of the classical Prandtl equation to viscoelastic fluids. It captures an interaction among viscous dissipation, convection and elasticity stretching
near the boundary. It can be formally derived from the boundary layer asymptotic analysis of viscoelastic fluids in the inviscid limit; see the Appendix. In this paper, we establish well-posedness of the system and investigate how the elasticity influences the stability of the boundary layer.

Prandtl's boundary layer theory has been a cornerstone of modern fluid mechanics and has a wide application across physical and engineering disciplines. The rigorous study of the Prandtl equation can be traced back to the pioneering work of Oleinik \cite{O63}, who established the local well-posedness for the 2D Prandtl equation under a monotone assumption on the horizontal velocity. Under the condition of a favorable pressure gradient, Xin and Zhang \cite{XZ04} constructed the global weak solutions to the Prandtl equations.  Subsequently, Alexander-Wang-Xu-Yang \cite{AWXY15} and Masmoudi-Wong \cite{MW15} independently developed the energy method to establish the local well-posedness in Sobolev spaces under the monotone condition.

When the monotone condition is violated, the Prandtl equation exhibits strong instabilities at high frequencies. In this direction, G\'erard-Varet and Dormy \cite{GD10} demonstrated Gevrey-2 ill-posedness for the linearized Prandtl equations around non-monotone shear flows. Further developments on nonlinear instability can be found in \cite{GN12, G00, GN11}.  Finite-time singularity is the other difficulty in the Prandtl equations.
 E and Engquist \cite{EE97} constructed a blow-up solution to the Prandtl equations. More recently, Kukavica-Vicol-Wang \cite{KVW17} investigated the Van Dommelen-Shen singularity.

To overcome these inherent instabilities,  an alternative approach is to establish well-posedness of the Prandtl equations in high-regularity frameworks, such as analytic or Gevrey spaces. In the pioneering work \cite{SC98}, Sammartino and Caflisch established local well-posedness for data analytic in tangential and normal variables. The analyticity requirement in the normal variable was later relaxed by Lombardo-Cannone-Sammartino \cite{LCS03}. The long-time well-posedness result was obtained by Ignotova-Vicol \cite{IV16} and Zhang-Zhang \cite{ZZ16}. More recently, Paicu-Zhang \cite{PZ21}  established the global-in-time well-posedness in the analytic framework. Beyond the analyticity setting, when the tangential velocity has a non-degenerate critical point, G\'erard-Varet and Masmoudi \cite{GM15} established well-posedness in Gevrey-$7/4$ space. This result was improved to the Gevrey-$2^-$ space by Chen-Wang-Zhang \cite{CWZ} and to Gevrey-$2$ by Li-Yang \cite{LY20}. Finally, Dietert-G\'erard-Varet \cite{DG19} established the local well-posedness in the critical Gevrey class 2 without any structural assumption. We also refer to \cite{CWZ22,GMM18,GMM24,KVW20,M14,NN18,SC98,SC982,WWZ17} for justification of the Prandtl ansatz and the inviscid limit of the Navier-Stokes equations.

The above results reveal that the stability of classical Prandtl boundary layer crucially depends on the monotonicity of the background velocity profiles and regularity of perturbations. A natural question is whether additional physical mechanisms can influence the stability property of the boundary layer. An important example is MHD boundary layer flow, where the non-degenerate tangential magnetic field can suppress the Prandtl-type instability. In this direction,
 G\'erard-Varet and Prestipino \cite{GP} derived several MHD boundary layer models and studied their linear stability. Liu-Xie-Yang \cite{LXY19,LXY20} established the well-posedness of MHD boundary layers and justified the corresponding Prandtl-type expansion. Furthermore, Liu-Wang-Xie-Yang \cite{LWXY20} studied the stability and instability mechanism of MHD boundary layers in the absence of magnetic resistivity.

In this paper, we aim to extend the mathematical theory of boundary layers to complex fluids. An important example is the viscoelastic model. Although the hydrodynamic stability in viscoelastic fluids has been extensively investigated in the physical literature, e.g., \cite{AH,R,TM}, a rigorous mathematical result on boundary layer theories has not been developed yet.  The coupling between fluid and elastic stress will introduce some new effects that are absent in the classical Newtonian setting. Understanding how these effects influence the stability of boundary layers is an important problem in mathematical theory of complex fluids.

We now introduce notation and norms in this paper. For a multi-index $\alpha=(\alpha_1, \alpha_2, \alpha_3)\in\mathbb{N}^3$, we define the differential operator:
\begin{align*}
D^{\alpha}=\partial_t^{\alpha_1}\partial_x^{\alpha_2}\partial_y^{\alpha_3},\qquad|\alpha|=\alpha_1+\alpha_2+\alpha_3.
\end{align*}
For a smooth function $f(t,x,y)$ and any integer $m\geq0$, we define the time-dependent energy norm:
\begin{align*}
\|f(t)\|_{H^m}=\left(\sum\limits_{|\alpha|\le m}\|D^{\alpha}f(t)\|_{L^2(\Omega)}^2\right)^{\frac12}.
\end{align*}
Note that the norm $\|\cdot\|_{H^m}$ includes both temporal and spatial derivatives. For a tangential multi-index $\alpha'=(\beta,\gamma)\in\mathbb{N}^2$, we use
\begin{align*}
\partial_\tau^{\alpha'}=\partial_t^{\beta}\partial_x^{\gamma},\quad |\alpha'|=\beta+\gamma
\end{align*}
to represent the tangential derivatives.  We define the mixed norm $\|\cdot\|_{p,q}$ as follows:
\begin{align*}
\|f(t)\|_{p, q}=\left(\sum\limits_{|\alpha|\le p \atop 0\le k\le q}\|\partial_\tau^{\alpha}\partial_y^kf(t, \cdot)\|_{L^2}^2\right)^{\frac12},
\end{align*}
where $p$ and $q$ are orders of tangential and normal derivatives, respectively. In particular, we use $\|f(t)\|_{m}$ to represent the norm $\|f(t)\|_{m, 0}$ which involves only tangential derivatives. It is straightforward to check that
\begin{align*}
\|f(t)\|_{H^m}^2\approx\sum\limits_{i=0}^m\|f(t)\|_{m-i, i}^2.
\end{align*}

\begin{theorem}\label{T1.1}
Let $m\geq 4$ be an integer. Let the initial data $(u_0, f_{1, 0}, f_{2, 0})$ satisfy the compatibility conditions up to order $m$,
\begin{align}
\left(u_0(x, y)-U(0, x), f_{1, 0}(x, y)-F_1(0, x), f_{2, 0}(x, y)-F_2(0, x)\right)\in H^{3m}(\Omega).\nonumber
\end{align}
and the initial deformation matrix is degenerate, namely,
\begin{equation}
\text{det}\left( \begin{array}{cc} f_{1,0}& f_{2,0}\\ f_{3,0}&f_{4,0} \end{array}\right)=0\quad\text{ in }\Omega,\label{d1}
\end{equation}
where initial data $f_{3,0}$ and $f_{4,0}$ are determined by the divergence free constraints in \eqref{1.1}. Assume further that the first row $(f_{1,0},f_{2,0})$ of the deformation matrix satisfies the following directional non-degeneracy: there exists a fixed non-zero vector $(k_1,k_2) \in \mathbb{R}^2$ and a constant $\sigma_0>0$ such that
\begin{align}\label{1.3}
k_1f_{1, 0}(x,y)+k_2f_{2, 0}(x,y)\ge 2\sigma_0,
\end{align}
for all $(x,y)\in \Omega$. There exists $T_*>0$, such that if prescribed outer flows $(U, P, F_1, F_2)(t, x)$ satisfy
\begin{align}
\sup\limits_{0\le t\le T_*}\sum\limits_{|\alpha|\le 2m}\|\partial_{\tau}^{\alpha}(U, P, F_1, F_2)(t, \cdot)\|_{L^2(\mathbb{T}_x)}<\infty,\nonumber
\end{align}
then the problem \eqref{1.1}-\eqref{1.2} admits a unique solution $(u, v, f_1, f_2, f_3, f_4)$ on the time interval $[0, T_*]$. Moreover, the solution satisfies
\begin{align}
\left(u(t, x, y)-U(t, x), f_{1}(t, x, y)-F_1(t, x), f_{2}(t, x, y)-F_2(t, x)\right)\in H^{m}(\Omega),\label{b1}
\end{align}
and
\begin{align}
k_1f_1(t, x, y)+k_2f_2(t, x, y)\ge \sigma_0, \label{1.3-1}
\end{align}
for all $(t, x, y)\in (0, T_*]\times\mathbb{T}\times\mathbb{R}_+$.
\end{theorem}

\begin{remark}
Our structural condition \eqref{1.3} only requires that the first row of the deformation tensor has a non-degenerate projection onto a fixed direction $(k_1,k_2)$. In particular, neither component is required to be non-degenerate, and either may vanish at some points. This contrasts with the MHD case studied in \cite{LWXY20,LXY19}, where the tangential magnetic field is required to be non-vanishing.
\end{remark}
\begin{remark}
The degeneracy \eqref{d1} on the deformation tensor is essential for the formation of the boundary layers and is preserved under the evolution; see Lemma \ref{lm2.1}.  In contrast, when the deformation tensor matrix is not singular near the boundary, the second and third authors \cite{WX23} proved that a boundary layer does not form and justified the inviscid limit of the compressible viscoelastic equations in strong norms.
\end{remark}

\begin{remark}
The rank-1 conditions \eqref{d1}-\eqref{1.3} are motivated by the strongly anisotropic stress in the high-Reynolds-number viscoelastic flows near solid walls. In this regime, the flow near the boundary is dominated by shear, and strong elastic effects can stretch the polymeric stress in the flow direction. Such background states have been extensively studied in the hydrodynamic stability theory of viscoelastic shear flows; see \cite{AH,R,TM} for their physical relevance.
\end{remark}

\begin{remark}
The initial data should satisfy the compatibility condition up to $m$-th order and required to lie in $H^{3m}$. Indeed, each time derivative can be converted into two spatial derivatives through the elasticity system. Hence, this requirement guarantees that the zeroth-order approximate solution constructed in Section $2.2$ belong to $H^m$.
\end{remark}

We now consider the case where both $f_1$ and $f_2$ vanish simultaneously at one point. In this case, the structural condition \eqref{1.3} is violated, and we expect instabilities to arise. To this end, we linearize the system \eqref{1.1} around a shear flow solution $(U_s(y), 0, F_{1, s}(y), F_{2, s}(y), 0, 0)$:
\begin{align}\label{1.5}
\begin{cases}
\partial_tu+U_s\partial_xu+v\partial_yU_s-\partial_y^2u-F_{1, s}\partial_xf_1-f_3F_{1, s}^{\prime}-F_{2, s}\partial_xf_2-f_4F_{2, s}^{\prime}=0,\\
\partial_tf_1+U_s\partial_xf_1+vF_{1, s}^{\prime}-F_{1, s}\partial_xu-f_3\partial_yU_s=0,\\
\partial_tf_2+U_s\partial_xf_2+vF_{2, s}^{\prime}-F_{2, s}\partial_xu-f_4\partial_yU_s=0,\\
\partial_xu+\partial_yv=0,\quad \partial_xf_1+\partial_yf_3=0,\quad \partial_xf_2+\partial_yf_4=0,\\
(u, v, f_3, f_4)|_{y=0}={\bf{0}},\quad \lim\limits_{y\to\infty}(u, f_1, f_2)(t, x, y)={\bf{0}}.
\end{cases}
\end{align}
%
For any $\alpha>0$, we introduce the weighted spaces
\begin{align*}
L^{ \infty}_\alpha(\mathbb{R}_+):= \{ g(y) \mid \|g\|_{L_{\alpha}^{\infty}}:=\sup_{y\geq  0}|e^{\alpha y}g(y)|<\infty \},
\end{align*}
and assume that $\partial_{y}^j\left(U_s(y)-U_\infty,F_{1,s}-F_{1,\infty},F_{2,s}-F_{2,\infty}\right)\in L^\infty_\alpha(\mathbb{R}_+)$, with $j=0,1,2,3.$ Here $U_\infty$, $F_{1,\infty}$ and $F_{2,\infty}$ are far fields of $U_s$, $F_{1,s}$ and $F_{2,s}$ respectively.
We define $$L^{2,\infty}_\alpha(\Omega):=L^2(\mathbb{T}, L^\infty_\alpha(\mathbb{R}_+))\subset L^2(\Omega).$$
For any $f\in L^{2}(\Omega)$, we denote by $\mathcal{P}_n (n\in \mathbb{Z})$  the orthogonal projection onto $n$-th Fourier mode in $x$:
$$
\left(\mathcal{P}_nf\right)(x,y)=f_{n}(y)e^{inx},~~\text{ where }f_n(y)=\frac{1}{2\pi}\int_0^{2\pi}f(x,y)e^{-inx}dx.
$$
We denote $\mathcal{S}(t)$ as the solution semigroup of the linearized equation \eqref{1.5}. Thanks to the shear flow structure of the background solution, we have $\mathcal{P}_n\mathcal{S}(t)=\mathcal{S}(t)\mathcal{P}_n$. We then define the projection $\mathcal{S}_n(t)$ onto the $n$-mode of $\mathcal{S}(t)$ as follows:
 $$\mathcal{S}_n(t):=\mathcal{S}(t)\mathcal{P}_n.$$

We have the following ill-posedness of the linearized problem \eqref{1.5}.
\begin{theorem}\label{T1.2}
Suppose that there exists a single point $a$, such that,
\begin{align*}
F_{1, s}(a)=F_{1, s}^{\prime}(a)=F_{2, s}(a)=F_{2, s}^{\prime}(a)=U_{s}^{\prime}(a)=0,\quad U_{s}^{\prime\prime}(a)\neq 0.
\end{align*}
Then there exists $\lambda>0$, such that for any $\mu\in [0,\frac14)$, and any $\delta>0,$ it holds that
\begin{align}\label{1.8}
\sup_{t\in [0,\delta],n\in \mathbb{Z}}\left(e^{-\lambda n^{\frac12}t}n^{-\mu}\|\mathcal{S}_{n}(t)\|_{L^{2,\infty}_\alpha(\Omega)\rightarrow L^{2,\infty}_\alpha(\Omega)}\right)=\infty.
\end{align}
\end{theorem}
\begin{remark}
Our results reveal the stabilizing effect of elasticity in the following aspects. On one hand, Theorem \ref{T1.1} establishes the stability of the boundary layer system under the directional non-degeneracy condition \eqref{1.3}, without any monotone condition on the velocity profile. On the other hand, if tangential components of the deformation matrix vanish simultaneously at one point, the elastic effects cannot suppress the Prandtl-type instability associated with a non-degenerate critical point of the velocity profile.
\end{remark}

We now explain the strategy of the proof. For the well-posedness theory, the main difficulty is that the system \eqref{1.1} exhibits a loss of derivatives. The most important step is therefore to eliminate these derivative-lossing terms by identifying a cancellation property. On the one hand, part of the derivative loss can be canceled through the coordinate transform \eqref{CT}, by using the stream function of the elastic variable together with the non-degeneracy condition on the first row of the deformation matrix \eqref{1.3}. On the other hand, owing to the degeneracy of the deformation matrix \eqref{d1}, the remaining derivative loss terms can also be canceled, so that the Sobolev solvability of solutions to \eqref{1.1} is established. This differs from \cite{LWXY20, LXY19}, where the derivative loss terms can be eliminated by only one mechanism.

In contrast, to establish the instability, we construct the approximate solution in the form \eqref{i9}, which gives rise to the Prandtl instability of the linearized system. It is noted that a key assumption is that both the shear flow $(F_{1, s}, F_{2, s})$ and its derivative vanish at a single point, so as to overcome the singularity in the error term of \eqref{i10}. From this viewpoint, the elasticity can be regarded as a perturbation of the velocity, and the instability mechanism originates from the Prandtl part. An interesting question is whether the requirement that the derivatives of
$(F_{1, s}, F_{2, s})$ also vanish at this single point is merely a technical restriction or is in fact essential.

The remainder of this paper is organized as follows. In Section \ref{S2}, we establish the local well-posedness of the viscoelastic boundary-layer model \eqref{1.2} in Sobolev spaces under the directional non-degeneracy condition. We first introduce an elastic stream function in Section \ref{S2.1}, and use it to reformulate the original system \eqref{1.2} as a hyperbolic-parabolic coupled system \eqref{EEE}.
We then derive the a priori bounds through successive $L^2$-estimates, tangential derivative estimates, and normal derivative estimates. Combining these estimates with the iteration scheme introduced in Section \ref{S2.2}, we obtain the local well-posedness result stated in Theorem  \ref{T1.1}. In Section \ref{S3}, we first construct an approximate growing mode, then establish the ill-posedness result in Theorem \ref{T1.2} by a contradiction argument. Finally, in the Appendix, we formally derive the viscoelastic Prandtl equations from a boundary layer expansion.

\bigskip

\section{Well-Posedness of Viscoelastic Boundary Layer Model} \label{S2}

\subsection{Coordinate transformation}\label{S2.1}
Let $(k_1,k_2) \in \mathbb{R}^2$ be the vector given in the Theorem \ref{T1.1}. From the divergence-free conditions in \eqref{1.1}, we have
\begin{align*}
\partial_x(k_1f_1+k_2f_2)+\partial_y(k_1f_3+k_2f_4)=0.
\end{align*}
Then there exists an elastic stream function $\phi$ satisfying
\begin{align}\label{2.1}
\partial_y\phi=k_1f_1+k_2f_2,\quad -\partial_x\phi=k_1f_3+k_2f_4.
\end{align}
Moreover, from the equations of $f_1$ and $f_2$ in \eqref{1.1} we deduce the following equation of $\phi$:
\begin{align}\label{2.2}
\partial_t\phi+u\partial_x\phi+v\partial_y\phi=0,
\end{align}
together with the Dirichlet boundary condition:
\begin{align}
\phi|_{y=0}=0.\nonumber
\end{align}
The initial data is given as follows:
\begin{align}
\phi_0(x, y):=\phi(0, x, y)=\int_0^y k_1f_{1, 0}(x, y^\prime)+k_2f_{2, 0}(x, y^\prime)\;dy^\prime.\nonumber
\end{align}
From the non-degenerate assumption \eqref{1.3}, we have
\begin{align}
\partial_y\phi_0(x, y)=k_1f_{1, 0}(x, y)+k_2f_{2, 0}(x, y)\ge 2\sigma_0>0.\nonumber
\end{align}
Hence $\phi_0(x, y)$ is a strictly increasing function with respect to $y$. Since $\phi(t, x, y)$ satisfies the transport equation \eqref{2.2},  it is also strictly increasing with respect to $y$.

It is well-known that the boundary layer model \eqref{1.1} exhibits a feature of loss of derivative in $x$ variable. To overcome this difficulty, we replace the normal variable $y$ with the stream function $\phi$, and introduce the following new coordinate system:
\begin{align}
\label{CT}
\bar{t}=t,\quad \bar{x}=x,\quad \bar{y}=\phi(t, x, y).
\end{align}
Under this change of coordinates, the space-time domain $[0,T]\times\Omega$ is mapped to itself, and the boundary $\{y=0\}$ is mapped to $\{\bar{y}=0\}$. Under this new coordinate system, the boundary layer model \eqref{1.1} can be reformulated as follows:
\begin{align}
\begin{cases}
\partial_{\bar{t}}u+u\partial_{\bar{x}}u+\cancel{(\partial_{t}\phi+u\partial_{x}\phi+v\partial_{y}\phi)}\partial_{\bar{y}}u+\partial_{\bar{x}}P-(k_1{f}_1+k_2{f}_2)\partial_{\bar{y}}((k_1{f}_1+k_2f_2)\partial_{\bar{y}}u)\\
\qquad-{f}_1\partial_{\bar{x}}{f}_1
-(f_1\partial_x\phi+f_3\partial_y\phi)\partial_{\bar{y}}{f}_1
-{f}_2\partial_{\bar{x}}{f}_2-(f_2\partial_x\phi+f_4\partial_y\phi)\partial_{\bar{y}}{f}_2=0,\\
\partial_{\bar{t}}{f}_1+{u}\partial_{\bar{x}}{f}_1+\cancel{(\partial_{t}\phi+u\partial_{x}\phi+v\partial_{y}\phi)}\partial_{\bar{y}}{f}_1
-{f}_1\partial_{\bar{x}}{u}
-(f_1\partial_x\phi+f_3\partial_y\phi)\partial_{\bar{y}}{u}=0,\\
\partial_{\bar{t}}{f}_2+{u}\partial_{\bar{x}}{f}_2+\cancel{(\partial_{t}\phi+u\partial_{x}\phi+v\partial_{y}\phi)}\partial_{\bar{y}}{f}_2
-{f}_2\partial_{\bar{x}}{u}
-(f_2\partial_x\phi+f_4\partial_y\phi)\partial_{\bar{y}}{u}=0.\nonumber
\end{cases}
\end{align}
Thanks to the equation \eqref{2.2}, the third terms (losing one derivative) in each of the above three equations vanish. Moreover, thanks to the non-degeneracy property \eqref{1.3-1}, the diffusion term in the equation for $u$ does not degenerate.

Now we study the last terms in each of the above three equations, which potentially lose one derivative. From \eqref{2.1}, we deduce that
\begin{align*}
f_1\partial_x\phi+f_3\partial_y\phi&=-k_2(f_1f_4-f_2f_3),\\
f_2\partial_x\phi+f_4\partial_y\phi&=k_1(f_1f_4-f_2f_3).
\end{align*}
\begin{lemma}[Propagation of determinant]\label{lm2.1}
Let  $(u,v,f_1,f_2,f_3,f_4)$ be a strong solution to \eqref{1.1}. Denote $\mathcal{D}:=f_1f_4-f_2f_3$ be the determinant of deformation matrix $\left( \begin{array}{cc} f_1& f_2\\ f_3&  f_4 \end{array}\right)$. Then it holds that:
\begin{align}
\partial_t\mathcal{D}+u\partial_x\mathcal{D}+v\partial_y\mathcal{D}=0.\label{z0}
\end{align}
In particular, if initially $\mathcal{D}_0=f_{1,0}f_{4,0}-f_{2,0}f_{3,0}\equiv0$, then $\mathcal{D}(t)\equiv 0$ for any $t>0.$
\end{lemma}
\begin{proof}
Using divergence-free conditions, we rewrite $\eqref{1.1}_2$ and $\eqref{1.1}_3$ as follows:
\begin{align}
&\partial_tf_1+\partial_y(vf_1-uf_3)=0,\label{z1}\\
&\partial_tf_2+\partial_y(vf_2-uf_4)=0.\label{z2}
\end{align}
Then from \eqref{z1} and \eqref{z2}, and using divergence free conditions for $(f_1,f_3)$ and $(f_2,f_4)$, we deduce the following evolution equations for $f_3$ and $f_4$:
\begin{align}
&\partial_tf_3+\partial_x(uf_3-vf_1)=0,\label{z3}\\
&\partial_tf_4+\partial_x(uf_4-uf_2)=0.\label{z4}
\end{align}

From \eqref{z1} and \eqref{z3}, we obtain
\begin{align}
\partial_t(f_1f_4)+u\partial_x(f_1f_4)+v\partial_y(f_1f_4)-f_3f_4\partial_yu-f_1f_2\partial_xv=0. \label{z6}
\end{align}
Similarly, from \eqref{z2} and \eqref{z4}, we obtain
\begin{align}
\partial_t(f_2f_3)+u\partial_x(f_2f_3)+v\partial_y(f_2f_3)-f_3f_4\partial_yu-f_1f_2\partial_xv=0. \label{z7}
\end{align}
Combining \eqref{z6} and \eqref{z7} together yields \eqref{z0}. The proof is complete.
\end{proof}

Thanks to Lemma \ref{lm2.1} and degenerate condition \eqref{d1},  all terms exhibiting loss of derivative vanish. Thus, we arrive at the following hyperbolic-parabolic coupled system:
\begin{align}
\label{EEE}
\begin{cases}
\partial_{\bar{t}}{u}+{u}\partial_{\bar{x}}{u}+\partial_{\bar{x}}p-(k_1{f}_1+k_2{f}_2)\partial_{\bar{y}}((k_1{f}_1+k_2{f}_2)\partial_{\bar{y}}{u})-{f}_1\partial_{\bar{x}}{f}_1
-{f}_2\partial_{\bar{x}}{f}_2=0,\\
\partial_{\bar{t}}{f}_1+{u}\partial_{\bar{x}}{f}_1
-{f}_1\partial_{\bar{x}}{u}=0,\\
\partial_{\bar{t}}{f}_2+{u}\partial_{\bar{x}}{f}_2
-{f}_2\partial_{\bar{x}}{u}=0,\\
{u}|_{\bar{y}=0}=0,\quad \lim\limits_{\bar{y}\to +\infty}({u}, {f}_1, {f}_2)=(U, F_1, F_2)(\bar{t}, \bar{x}),\\
({u}, {f}_1, {f}_2)(0, \bar{x}, \bar{y})=(u_0, f_{1, 0}, f_{2, 0})(\bar{x}, \bar{y}).
\end{cases}
\end{align}
Note that the assumption \eqref{1.3} not only ensures that the coordinate transformation \eqref{CT} is invertible, but also guarantees the positivity of the diffusion coefficient in \eqref{EEE}. Henceforth, for notational simplicity, we drop the bars and continue to denote the transformed coordinates by $(t,x,y)$ when no confusion arises.

To eliminate far-fields  $(U,F_1,F_2)$, we define a cut-off function $\chi({y})\in C^{\infty}(\mathbb{R}_+)$ satisfying
\begin{align}\nonumber
\chi({y})=
\begin{cases}
1,\quad \mbox{if}\,\, {y}\ge 2,\\
0,\quad \mbox{if}\,\, 0\le {y}\le 1,
\end{cases}
\end{align}
and then introduce the new unknowns
\begin{align*}
\hat{u}(t, x, y)&={u}(t, x, y)-U(t, x)\chi(y),\\
\hat{f}_1(t, x, y)&={f}_1(t, x, y)-F_1(t, x),\\
\hat{f}_2(t, x, y)&={f}_2(t, x, y)-F_2(t, x).
\end{align*}
Then $(\hat{u},\hat{f}_1,\hat{f}_2)$ satisfies the following system:
\begin{align}\label{2.9}
\begin{small}
\begin{cases}
\partial_t\hat{u}+(\hat{u}+\chi U)\partial_x\hat{u}-(\hat{f}_1+F_1)\partial_x\hat{f}_1-(\hat{f}_2+F_2)\partial_x\hat{f}_2+\hat{u}\chi\partial_xU -\hat{f}_1(\partial_xF_1+\chi^{\prime\prime}UF_1)-\hat{f}_2\partial_xF_2\\
\qquad\qquad-\left[(k_1(\hat{f}_1+F_1)+k_2(\hat{f}_2+F_2)\right]\partial_y\left[\left(k_1(\hat{f}_1+F_1)+k_2(\hat{f}_2+F_2)\right)\partial_y\hat{u}+\chi^
\prime U\hat{f}_1\right]
=r_1,\\
\partial_t\hat{f}_1+(u+\chi U)\partial_x\hat{f}_1-(\hat{f}_1+F_1)\partial_x\hat{u}+\hat{u}\partial_xF_1-\hat{f}_1\chi\partial_xU=r_2,\\
\partial_t\hat{f}_2+(u+\chi U)\partial_x\hat{f}_2-(\hat{f}_2+F_2)\partial_x\hat{u}+\hat{u}\partial_xF_2-\hat{f}_2\chi\partial_xU=r_3,\\
\hat{u}|_{y=0}=0,\quad \lim\limits_{y\to\infty}(\hat{u}, \hat{f}_1, \hat{f}_2)(t, x, y)=(0, 0, 0),\\
(\hat{u}, \hat{f}_1, \hat{f}_2)|_{t=0}=(u_0(x, y)-\chi(y)U(0, x), f_{1, 0}(x, y)-F_1(0, x), f_{2, 0}(x, y)-F_2(0, x)),
\end{cases}
\end{small}
\end{align}
where the error terms $(r_1,r_2,r_3)$ are defined by
\begin{align*}
r_1&=-\chi\partial_tU-\chi^2U\partial_xU+F_1\partial_xF_1+F_2\partial_xF_2+\chi^{\prime\prime}UF_1^2-\partial_xP\\
&=(1-\chi)\partial_tU+(1-\chi^2)U\partial_xU+\chi^{\prime\prime}UF_1^2,\\
r_2&=-\partial_tF_1-\chi U\partial_xF_1+\chi F_1\partial_xU\\
&=-(1-\chi)\partial_tF_1,\\
r_3&=-\partial_tF_2-\chi U\partial_xF_2+\chi F_2\partial_xU\\
&=-(1-\chi)\partial_tF_2,
\end{align*}
where the Bernoulli law \eqref{1.2} has been used. Henceforth, for notational simplicity, we drop the hats and continue to denote the perturbation variables by $(u,f_1,f_2)$ when no confusion arises.

Now we establish the solvability of \eqref{2.9}.

\begin{proposition}\label{T2.1}
Let $m\ge 4$ be an integer and $(k_1,k_2)\in \mathbb{R}^2$ be the vector given in Theorem \ref{T1.1}. Suppose the initial data $(u_0, f_{1, 0}, f_{2, 0})$ satisfy compatibility conditions up to $m$-th order,
\begin{align}
(u_0(x, y)-U(0, x), f_{1, 0}(x, y)-F_1(0, x), f_{2, 0}(x, y)-F_2(0, x))\in H^{3m}(\Omega),\nonumber
\end{align}
and there exists a positive constant $\sigma_0$, such that
\begin{align*}
k_1(f_{1, 0}(x, y)+F_1(0, x))+k_2(f_{2, 0}(x, y)+F_2(0, x))\geq 2\sigma_0>0.
\end{align*}
Then there exists $T_*>0$, such that if the outflow $(U, P, F_1, F_2)(t, x)$ satisfies
\begin{align}\label{2.10}
\sup\limits_{0\le t\le T_*}\sum\limits_{|\alpha|\le 2m}\|\partial_{\tau}^{\alpha}(U, P, F_1, F_2)(t, \cdot)\|_{L^2(\mathbb{T}_x)}<\infty,
\end{align}
 the system \eqref{2.9} admits a unique solution $(u, f_1, f_2)$ on $[0, T_*]$. Moreover, the solution satisfies
\begin{align}
\left(u(t, x, y)-U(t, x), f_{1}(t, x, y)-F_1(t, x), f_{2}(t, x, y)-F_2(t, x)\right)\in H^{m}(\Omega),\label{b}
\end{align}
for any $t>0$, and
\begin{align}\label{2.13}
f_k\eqdef k_1(f_1+F_1)+(k_2f_2+F_2)\ge \sigma_0.
\end{align}
\end{proposition}
The proof of Proposition \ref{T2.1} relies on the Picard iteration and the fixed point theorem. We divide it into the following subsections.

\subsection{Iteration scheme}\label{S2.2}

We construct the approximate solutions by an iteration scheme.

\noindent
 \underline{Initial Step:}
 We first construct the zeroth-order approximate solution that satisfies $m$-th order compatibility conditions and the non-degenerate condition \eqref{2.13}. For this, we denote
\begin{align*}
(u_0^j, f_{1, 0}^j, f_{2, 0}^j)(x, y)=(\partial_t^ju, \partial_t^jf_1, \partial_t^jf_2)(0, x, y),\quad 0\le j\le m,
\end{align*}
where the time derivatives are inductively obtained from the equation \eqref{2.9}.  We then construct the zeroth-order approximate solution as follows:
\begin{align}
u^0(t, x, y)=\sum\limits_{j=0}^{m}\frac{t^j}{j!}u_0^j(x, y),\quad f_1^0(t, x, y)=\sum\limits_{j=0}^{m}\frac{t^j}{j!}f_{1, 0}^j(x, y),\quad f_2^0(t, x, y)=\sum\limits_{j=0}^{m}\frac{t^j}{j!}f_{2, 0}^j(x, y).\nonumber
\end{align}
It is straightforward to check that there exists a time $T_0\in (0, T]$ such that for any $t\in (0,T_0)$, the solution $(u^0,f_1^0,f^2_0)$ satisfies \eqref{b}, and
\begin{align*}
f_k^0\eqdef k_1(f_1^0(t, x, y)+F_1(t, x))+k_2(f_2^0(t, x, y)+F_2(t, x))\ge \sigma_0.
\end{align*}

\noindent
\underline{Iteration:} 
Next, we construct the $(n+1)$-th order approximate solution $(u^{n+1}, f_1^{n+1}, f_2^{n+1})$ by induction. Assume that the $n$-th order approximate solution $(u^{n}, f_1^{n}, f_2^{n})$ is defined over the time interval $(0, T_n]$ for some $T_n\le T_0$, and satisfies
\begin{align}\label{2.15}
f_k^n\eqdef k_1(f_1^n(t, x, y)+F_1(t, x))+k_2(f_2^n(t, x, y)+F_2(t, x))\ge \sigma_0.
\end{align}
Then we define $(n+1)$-th order approximate solution $(u^{n+1}, f_1^{n+1}, f_2^{n+1})$ as the solution to the following linearized system:
\begin{align}\label{2.16}
\begin{cases}
\partial_tu^{n+1}+(u^n+\chi U)\partial_xu^{n+1}-f_k^n\partial_y(f_k^n\partial_yu^{n+1}+\chi^\prime Uf_1^{n+1})-(f_1^n+F_1)\partial_xf_1^{n+1}\\
\quad\qquad-(f_2^n+F_2)\partial_xf_2^{n+1}+\chi\partial_xU u^{n+1}-(\partial_xF_1+\chi^{\prime\prime}UF_1)f_1^{n+1}-\partial_xF_2f_2^{n+1}=r_1,\\
\partial_tf_1^{n+1}+(u^n+\chi U)\partial_xf_1^{n+1}-(f_1^n+F_1)\partial_xu^{n+1}+\partial_xF_1u^{n+1}-\chi\partial_xUf_1^{n+1}=r_2,\\
\partial_tf_2^{n+1}+(u^n+\chi U)\partial_xf_2^{n+1}-(f_2^n+F_2)\partial_xu^{n+1}+\partial_xF_2u^{n+1}-\chi\partial_xUf_2^{n+1}=r_3,\\
u^{n+1}|_{y=0}=0,\quad \lim\limits_{y\to\infty}(u^{n+1}, f_1^{n+1}, f_2^{n+1})(t, x, y)=(0, 0, 0),\\
(u^{n+1}, f_1^{n+1}, f_2^{n+1})|_{t=0}=(u_0(x, y)-\chi(y)U(0, x), f_{1, 0}(x, y)-F_1(0, x), f_{2, 0}(x, y)-F_2(0, x)).
\end{cases}
\end{align}
Thanks to  \eqref{2.10}, the error terms on the right hand side of \eqref{2.16} satisfy
\begin{align}\label{2.17}
\|(r_1, r_2, r_3)\|_{H^m}\le C.
\end{align}

\subsection{$L^2$ estimates}
Let ${\bf \Lambda}=(u,f_1,f_2)$, and ${\bf \Lambda}^n=(u^n,f_1^n,f_2^n)$ with $n=1,2,\cdots$.
\begin{lemma}\label{P2.2}
There exists a positive constant $C$ depending on $\sigma_0$, such that the solution to \eqref{2.16} satisfies
\begin{align}
\frac{d}{dt}\|{\bf \Lambda}^{n+1}\|_{L^2}^2+\|\partial_yu^{n+1}\|_{L^2}^2\le C(1+\|{\bf \Lambda}^n\|_{W^{1, \infty}}^2)\|{\bf \Lambda}^{n+1}\|_{L^2}^2+C.\label{ad1}
\end{align}
\end{lemma}

\begin{proof}
Taking inner product of \eqref{2.16} with ${\bf \Lambda}^{n+1}=(u^{n+1}, f_1^{n+1}, f_2^{n+1})$, we obtain
\begin{align}
&\frac{1}{2}\frac{d}{dt}\|(u^{n+1}, f_1^{n+1}, f_2^{n+1})\|_{L^2}^2+\sum_{i=1}^5 I_i=0,\label{a2.1}
\end{align}
where
\begin{align}
I_1=&\langle (u^n+\chi U)\partial_xu^{n+1}, u^{n+1}\rangle+\langle (u^n+\chi U)\partial_xf_1^{n+1}, f_1^{n+1}\rangle+\langle (u^{n}+\chi U)\partial_xf_2^{n+1}, f_2^{n+1}\rangle,\notag\\
I_2=&-\langle(f_1^{n}+F_1)\partial_xf_1^{n+1},u^{n+1}\rangle+\langle(f_2^{n}+F_2)\partial_xf_2^{n+1}, u^{n+1}\rangle\nonumber\\
&-\langle (f_1^{n}+F_1)\partial_xu^{n+1}, f_1^{n+1}\rangle-\langle (f_2^{n}+F_2)\partial_xu^{n+1}, f_2^{n+1}\rangle,\nonumber\\
I_3=&\langle \chi\partial_xU u^{n+1}-(\partial_xF_1+\chi^{\prime\prime}UF_1)f_1^{n+1}-\partial_xF_2f_2^{n+1}, u^{n+1}\rangle\notag\\
&+\langle \partial_xF_1u^{n+1}-\chi\partial_xUf_1^{n+1}, f_1^{n+1}\rangle+\langle\partial_xF_2u^{n+1}-\chi\partial_xUf_2^{n+1}, f_2^{n+1}\rangle,\notag\\
I_4=&-\langle f_k^{n}\partial_y(f_k^{n}\partial_yu^{n+1}+\chi^\prime Uf_1^{n+1}), u^{n+1}\rangle,\notag\\
I_5=&-\langle r_1, u^{n+1}\rangle-\langle r_2, f_1^{n+1}\rangle-\langle r_3, f_2^{n+1}\rangle.\nonumber
\end{align}

Now we estimate $I_1$--$I_5$ term by term. Integrating by parts, we obtain
\begin{align}\label{2.6}
|I_1|&=\frac{1}{2}\left|-\int_{\Omega}\partial_x(u^n+\chi U)(|u^{n+1}|^2+|f_1^{n+1}|^2+|f_2^{n+1}|^2)\;dxdy\right|\notag\\
&\le C(1+\|\partial_xu^n\|_{L^\infty})\|{\bf\Lambda}^{n+1}\|_{L^2}^2.
\end{align}
Similarly, for $I_2$, we integrate by parts to obtain that
\begin{align}\label{2.7-1}
|I_2|&=\left|\int_{\Omega}\partial_x(f_1^n+F_1)u^{n+1}f_1^{n+1}+\partial_x(f_2^n+F_2)u^{n+1}f_2^{n+1}\;dxdy \right|\notag\\
&\le C(1+\|(\partial_xf_1^n, \partial_xf_2^n)\|_{L^\infty})\|{\bf \Lambda}^{n+1}\|_{L^2}^2.
\end{align}

For $I_3$ involving the zeroth-order term and $I_5$ involving error terms, we deduce that
\begin{align}\label{2.8}
|I_3|+|I_5|\le C\|{\bf \Lambda}^{n+1}\|_{L^2}^2+C\|(r_1, r_2, r_3)\|_{L^2}\|{\bf\Lambda}^{n+1}\|_{L^2}\leq C\|{\bf\Lambda}^{n+1}\|_{L^2}^2+C.
\end{align}

For the diffusion term $I_4$, we integrate by parts and use the no-slip boundary condition $u^{n+1}|_{y=0}=0$ as well as non-degenerate condition \eqref{2.15} to obtain
\begin{align}
I_4&=\|f_k^n\partial_yu^{n+1}\|_{L^2}^2+\int_{\Omega}\partial_yf_k^nf_k^n\partial_yu^{n+1}u^{n+1}\;dxdy\notag\\
&\quad+\int_{\Omega}\partial_yf_k^n\chi^{\prime}Uf_1^{n+1}u^{n+1}\;dxdy+\int_{\Omega}f_k^n\chi^{\prime}Uf_1^{n+1}\partial_yu^{n+1}\;dxdy\notag\\
&\ge \frac{1}{2}\|f_k^n\partial_yu^{n+1}\|_{L^2}^2-\frac{3}{2}\int_{\Omega}(\partial_yf_k^nu^{n+1})^2+(\chi^{\prime}Uf_1^{n+1})^2\;dxdy\notag\\
&\ge \frac{\sigma_0}{2}\|\partial_yu^{n+1}\|_{L^2}^2-C(1+\|(\partial_yf_1^n, \partial_yf_2^n)\|_{L^\infty}^2)\|\Lambda^{n+1}\|_{L^2}^2.\label{2.8-1}
\end{align}
Substituting \eqref{2.6}-\eqref{2.8-1} into \eqref{a2.1} yields \eqref{ad1}. The proof of Lemma \ref{P2.2} is complete.
\end{proof}

\subsection{Tangential estimates}
\begin{lemma}\label{P2.3}
The solution to \eqref{2.16} satisfies
\begin{align}
\frac{d}{dt}\|{\bf\Lambda}^{n+1}\|_{m}^2+\|\partial_yu^{n+1}\|_{m}^2\le C(1+\|{\bf\Lambda}^n\|_{H^m}^4)\|{\bf\Lambda}^{n+1}\|_{H^m}^2+C.\label{ad2}
\end{align}
\end{lemma}

\begin{proof}
Applying the tangential derivative $\partial_\tau^\alpha$ to the first three equations in \eqref{2.16} with $|\alpha|\le m$, then taking inner product of the resulting equations with $(\partial_\tau^{\alpha}u^{n+1}, \partial_\tau^{\alpha}f_1^{n+1}, \partial_\tau^{\alpha}f_2^{n+1})$, we obtain the following identity:
\begin{align}
\label{EQ}
\frac{1}{2}\frac{d}{dt}\|(\partial_\tau^{\alpha}u^{n+1}, \partial_\tau^{\alpha}f_1^{n+1}, \partial_\tau^{\alpha}f_2^{n+1})\|_{L^2}^2+\sum_{i=1}^8J_i=0,
\end{align}
where $J_1$--$J_8$ are defined as follows:
\begin{align}
J_1&=\langle (u^n+\chi U)\partial_x\partial_\tau^{\alpha}u^{n+1}, \partial_\tau^{\alpha}u^{n+1}\rangle+\langle (u^n+\chi U)\partial_x\partial_\tau^{\alpha}f_1^{n+1}, \partial_\tau^{\alpha}f_1^{n+1}\rangle\nonumber\\
&\quad+\langle (u^{n}+\chi U)\partial_x\partial_\tau^{\alpha}f_2^{n+1}, \partial_\tau^{\alpha}f_2^{n+1}\rangle,\notag\\
J_2&=\langle\mathcal{C}_1\partial_xu^{n+1}, \partial_\tau^{\alpha}u^{n+1}\rangle+\langle\mathcal{C}_1\partial_xf_1^{n+1}, \partial_\tau^{\alpha}f_1^{n+1}\rangle+\langle\mathcal{C}_1\partial_xf_2^{n+1}, \partial_\tau^{\alpha}f_2^{n+1}\rangle,\notag\\
J_3&=-\langle (f_1^{n}+F_1)\partial_x\partial_\tau^{\alpha}f_1^{n+1}, \partial_\tau^{\alpha}u^{n+1}\rangle-\langle (f_1^{n}+F_1)\partial_x\partial_\tau^{\alpha}u^{n+1}, \partial_\tau^{\alpha}f_1^{n+1}\rangle
\notag\\
&\quad-\langle(f_2^{n}+F_2)\partial_x\partial_\tau^{\alpha}f_2^{n+1}, \partial_\tau^{\alpha}u^{n+1}\rangle-\langle (f_2^{n}+F_2)\partial_x\partial_\tau^{\alpha}u^{n+1}, \partial_\tau^{\alpha}f_2^{n+1}\rangle,\notag\\
J_4&=-\langle \mathcal{C}_2\partial_xf_1^{n+1}, \partial_\tau^{\alpha}u^{n+1}\rangle- \langle \mathcal{C}_2\partial_xu^{n+1}, \partial_\tau^{\alpha}f_1^{n+1}\rangle\nonumber\\
&\quad-\langle \mathcal{C}_3\partial_xf_2^{n+1}, \partial_\tau^{\alpha}u^{n+1}\rangle-\langle \mathcal{C}_3\partial_xu^{n+1}, \partial_\tau^{\alpha}f_2^{n+1}\rangle,\notag\\
J_5&=\langle \partial_\tau^{\alpha}(\chi\partial_xU u^{n+1}-(\partial_xF_1+\chi^{\prime\prime}UF_1)f_1^{n+1}-\partial_xF_2f_2^{n+1}), \partial_\tau^{\alpha}u^{n+1}\rangle\notag\\
&\quad+\langle \partial_\tau^{\alpha}(\partial_xF_1u^{n+1}-\chi\partial_xUf_1^{n+1}), \partial_\tau^{\alpha}f_1^{n+1}\rangle+\langle \partial_\tau^{\alpha}(\partial_xF_2u^{n+1}-\chi\partial_xUf_2^{n+1}), \partial_\tau^{\alpha}f_2^{n+1}\rangle,\notag\\
J_6&=-\langle f_k^{n}\partial_y\partial_\tau^{\alpha}(f_k^{n}\partial_yu^{n+1}+\chi^\prime Uf_1^{n+1}), \partial_\tau^{\alpha}u^{n+1}\rangle,\notag\\
J_7&=-\langle \mathcal{C}_k\partial_y(f_k^{n}\partial_yu^{n+1}+\chi^\prime Uf_1^{n+1}), \partial_\tau^{\alpha}u^{n+1}\rangle,\notag\\
J_8&=-\langle \partial_\tau^{\alpha}r_1, \partial_\tau^{\alpha}u^{n+1}\rangle-\langle \partial_\tau^{\alpha}r_2, \partial_\tau^{\alpha}f_1^{n+1}\rangle-\langle \partial_\tau^{\alpha}r_3, \partial_\tau^{\alpha}f_2^{n+1}\rangle.\nonumber
\end{align}
Here $\mathcal{C}_{i}(i=1, 2, 3)$ denotes the commutators
\begin{align}
\mathcal{C}_1=[\partial_\tau^{\alpha}, u^n+\chi U],\quad \mathcal{C}_2=[\partial_\tau^{\alpha}, f_1^n+F_1],\quad \mathcal{C}_3=[\partial_\tau^{\alpha}, f_2^n+F_2],\nonumber
\end{align}
and $\mathcal{C}_{k}$ is the linear combination of $\mathcal{C}_2, \mathcal{C}_3$
\begin{align}
\mathcal{C}_k=k_1\mathcal{C}_2+k_2\mathcal{C}_3,\nonumber
\end{align}
where $(k_1,k_2)$ is the vector given in Theorem \ref{T1.1}.

Now we begin to estimate $J_1$--$J_8$. Integrating by parts yields
\begin{align}
|J_1|&=\frac{1}{2}\left|-\int_{\Omega}\partial_x(u^n+\chi U)(|\partial_\tau^{\alpha}u^{n+1}|^2+|\partial_\tau^{\alpha}f_1^{n+1}|^2+|\partial_\tau^{\alpha}f_2^{n+1}|^2)\;dxdy\right|\notag\\
&\le C(1+\|\partial_xu^n\|_{L^\infty})\|\partial_\tau^{\alpha}{\bf\Lambda}^{n+1}\|_{L^2}^2.\label{P2.3-1}
\end{align}
Next we estimate $J_2$ and $J_4$ involving commutators. For any function $f$,  we compute
\begin{align}
\|\mathcal{C}_1f\|_{L^2}&=\|[\partial_\tau^{\alpha}, u^n+\chi U]f\|_{L^2}\leq C\sum\limits_{\beta+\gamma=\alpha \atop \beta\neq 0}\|\partial_\tau^{\beta}(u^n+\chi U)\partial_\tau^{\gamma}f\|_{L^2}\nonumber\\
&\leq C(1+\|u^n\|_{H^m})\|f\|_{H^{m-1}}.\label{com1}
\end{align}
Similarly, we obtain
\begin{align}
\|\mathcal{C}_2f\|_{L^2}+\|\mathcal{C}_3f\|_{L^2}\leq C(1+\|f_1^n\|_{H^m}+\|f_2^n\|_{H^m})\|f\|_{H^{m-1}}.\label{com}
\end{align}

Substituting \eqref{com1} and \eqref{com} into $J_2$ and $J_4$, we obtain
\begin{align}\label{P2.3-2}
|J_2|+|J_4|&\leq C\left(\sum_{i=1}^3\|\mathcal{C}_i\partial_x{\bf\Lambda}^{n+1}\|_{L^2}\right)\|\partial_\tau^{\alpha}{\bf\Lambda}^{n+1}\|_{L^2}\le C(1+\|{\bf \Lambda}^n\|_{H^m})\|{\bf\Lambda}^{n+1}\|_{m}^2.
\end{align}
For $J_3$, we integrate by parts to obtain
\begin{align}
|J_3|&=\left|\int_{\Omega}\partial_x(f_1^n+F_1)|\partial_\tau^{\alpha}f_1^{n+1}|^2\;dxdy\right|+\left|\int_{\Omega}\partial_x(f_2^n+F_2)|\partial_\tau^{\alpha}f_2^{n+1}|^2\;dxdy\right|\notag\\
&\le C(1+\|(\partial_xf_1^n, \partial_xf_2^n)\|_{L^\infty})\|\partial_\tau^{\alpha}{\bf\Lambda}^{n+1}\|_{L^2}^2.\label{P2.3-3}
\end{align}
By a direct computation, we can control $J_5$ and $J_8$, which respectively involve lower-order terms and error terms as follows:
\begin{align}
|J_5|+|J_8|\le &C\|{\bf\Lambda}^{n+1}\|_{m}^2+C\|(r_1, r_2,r_3)\|_{m}\|{\bf\Lambda}^{n+1}\|_{m}\nonumber\\
\leq & C\|{\bf\Lambda}^{n+1}\|_{m}^2+C.\label{P2.3-4}
\end{align}

For the diffusion term $J_6$, we integrate by parts and use the boundary condition $\partial_\tau^{\alpha}u^{n+1}|_{y=0}=0$ to get
\begin{align*}
J_6&=\langle f_k^n\partial_\tau^{\alpha}(f_k^n\partial_yu^{n+1}+\chi^\prime Uf_1^{n+1}), \partial_\tau^{\alpha}\partial_yu^{n+1}\rangle\\
&\quad+\langle \partial_yf_k^n\partial_\tau^{\alpha}(f_k^n\partial_yu^{n+1}+\chi^\prime Uf_1^{n+1}), \partial_\tau^{\alpha}u^{n+1}\rangle\\
&=: J_{61}+J_{62}.
\end{align*}
For the first term, note that
\begin{align*}
\partial_\tau^{\alpha}(f_k^n\partial_yu^{n+1})=f_k^n\partial_\tau^{\alpha}\partial_yu^{n+1}+\mathcal{C}_k\partial_yu^{n+1}.
\end{align*}
Then using the non-degeneracy \eqref{2.15} and commutator estimates \eqref{com}, we obtain
\begin{align}\label{2.40}
J_{61}&=\|f_k^n\partial_\tau^{\alpha}\partial_yu^{n+1}\|_{L^2}^2+\int_{\Omega}f_k^n(\mathcal{C}_k\partial_yu^{n+1}+\partial_\tau^{\alpha}(\chi^{\prime}Uf_1^{n+1}))\partial_\tau^{\alpha}\partial_yu^{n+1}\;dxdy\notag\\
&\ge \frac{3}{4}\|f_k^n\partial_\tau^{\alpha}\partial_yu^{n+1}\|_{L^2}^2-\int_{\Omega}(\mathcal{C}_k\partial_yu^{n+1})^2+(\partial_\tau^{\alpha}(\chi^{\prime}Uf_1^{n+1}))^2\;dxdy\notag\\
&\ge \frac{3}{4}\|f_k^n\partial_\tau^{\alpha}\partial_yu^{n+1}\|_{L^2}^2-C(1+\|(f_1^n, f_2^n)\|_{H^m}^2)(\|f_1^{n+1}\|_{H^m}^2+\|\partial_yu^{n+1}\|_{H^{m-1}}^2).
\end{align}
For $J_{62}$, we have
\begin{align*}
J_{62}&=\int_{\Omega}\partial_yf_k^nf_k^n\partial_\tau^{\alpha}\partial_yu^{n+1}\partial_\tau^{\alpha}u^{n+1}\;dxdy+\int_{\Omega}\partial_yf_k^n\mathcal{C}_k\partial_yu^{n+1}\partial_\tau^{\alpha}u^{n+1}\;dxdy\notag\\
&\quad+\int_{\Omega}\partial_yf_k^n\partial_\tau^{\alpha}(\chi^{\prime}Uf_1^{n+1})\partial_\tau^{\alpha}u^{n+1}\;dxdy\notag\\
&\le \frac{1}{4}\|f_k^n\partial_\tau^{\alpha}\partial_yu^{n+1}\|_{L^2}^2+C\int_{\Omega}(\partial_yf_k^n\partial_\tau^{\alpha}u^{n+1})^2\;dxdy\nonumber\\
&\quad+C\int_{\Omega}(\mathcal{C}_k\partial_yu^{n+1})^2+(\partial_\tau^{\alpha}(\chi^{\prime}Uf_1^{n+1}))^2\;dxdy\nonumber\\
&\leq \frac{1}{4}\|f_k^n\partial_\tau^{\alpha}\partial_yu^{n+1}\|_{L^2}^2+C(1+\|(f_1^n, f_2^n)\|_{H^m}^2)(\|f_1^{n+1}\|_{H^m}^2+\|\partial_yu^{n+1}\|_{H^{m-1}}^2).
\end{align*}
This together with \eqref{2.40} implies
\begin{align}
J_{6}\ge \frac{1}{2}\|f_k^n\partial_\tau^{\alpha}\partial_yu^{n+1}\|_{L^2}^2-C(1+\|(f_1^n, f_2^n)\|_{H^m}^2)(\|(f_1^{n+1}, u^{n+1})\|_{H^m}^2+\|\partial_yu^{n+1}\|_{H^{m-1}}^2).\label{P2.3-5}
\end{align}

Finally, from the definition of the commutator, we can separate $J_7$ into two parts.
\begin{align*}
J_7&=\sum\limits_{\beta+\gamma=\alpha \atop |\beta|=1}\langle m\partial_\tau^{\beta}f_k^{n}\partial_\tau^{\gamma}\partial_y(f_k^n\partial_yu^{n+1}+\chi^\prime Uf_1^{n+1}), \partial_\tau^{\alpha}u^{n+1}\rangle\\
&\quad+\sum\limits_{\beta+\gamma=\alpha \atop |\beta|\ge 2}\langle c_{\beta, \gamma}\partial_\tau^{\beta}f_k^{n}\partial_\tau^{\gamma}\partial_y(f_k^n\partial_yu^{n+1}+\chi^\prime Uf_1^{n+1}), \partial_\tau^{\alpha}u^{n+1}\rangle\\
&=:J_{71}+J_{72}.
\end{align*}
Integrating by parts and using the boundary condition $\partial_\tau^{\alpha}u^{n+1}|_{y=0}=0$ again, we obtain
\begin{align*}
J_{71}&=-\sum\limits_{\beta+\gamma=\alpha \atop |\beta|=1}\langle m\partial_\tau^{\beta}\partial_yf_k^{n}\partial_\tau^{\gamma}(f_k^n\partial_yu^{n+1}+\chi^\prime Uf_1^{n+1}), \partial_\tau^{\alpha}u^{n+1}\rangle\\
&\quad -\sum\limits_{\beta+\gamma=\alpha \atop |\beta|=1}\langle m\partial_\tau^{\beta}f_k^{n}\partial_\tau^{\gamma}(f_k^n\partial_yu^{n+1}+\chi^\prime Uf_1^{n+1}), \partial_\tau^{\alpha}\partial_yu^{n+1}\rangle.
\end{align*}
Then a direct computation yields that
\begin{align}
|J_{71}|+|J_{72}|&\le C(1+\|(f_1^n, f_2^n)\|_{H^m}^2)(\|\partial_yu^{n+1}\|_{H^{m-1}}+\|f_1^{n+1}\|_{H^{m-1}})\|u^{n+1}\|_{m}\notag\\
&\quad+C(1+\|(f_1^n, f_2^n)\|_{H^m}^4)\|{\bf\Lambda}^{n+1}\|_{H^m}^2+\frac{1}{4}\|\partial_\tau^{\alpha}\partial_yu^{n+1}\|_{L^2}\nonumber\\
&\leq \frac{1}{4}\|\partial_\tau^{\alpha}\partial_yu^{n+1}\|_{L^2}+C(1+\|(f_1^n, f_2^n)\|_{H^m}^4)\|{\bf\Lambda}^{n+1}\|_{H^m}^2.\label{P2.3-6}
\end{align}
Substituting \eqref{P2.3-1}-\eqref{P2.3-6} into \eqref{EQ}, we obtain \eqref{ad2}. The proof of Lemma \ref{P2.3} is complete.
\end{proof}

\subsection{Estimates of first-order normal derivative of velocity}
\begin{lemma}\label{P2.4}
There exists a positive constant $\delta$, which can be arbitrarily small, such that for any solutions ${\bf \Lambda}^{n+1}$ to \eqref{2.16}, it holds
\begin{align}
&\frac{d}{dt}\sum\limits_{|\alpha|\le m-1}\|f_k^n\partial_\tau^{\alpha}\partial_yu^{n+1}\|_{L^2}^2+\|\partial_tu^{n+1}\|_{m-1}^2\notag\\
&\quad \le  \delta\|\partial_yu^{n+1}\|_{m}^2+C_\delta(1+{\bf \|\Lambda}^n\|_{H^{m-1}}^4)\|{\bf \Lambda}^{n+1}\|_{H^{m}}^2+C_\delta.\label{ad3}
\end{align}
\end{lemma}

\begin{proof}
Applying the operator $\partial_\tau^\alpha$ to the first equation in \eqref{2.16} with $|\alpha|\le m-1$, then taking inner product of the resulting equations with $\partial_t\partial_\tau^{\alpha}u^{n+1}$, we obtain that
\begin{align}
&\|\partial_\tau^\alpha\partial_tu^{n+1}\|_{L^2}^2+
\underbrace{\langle \partial_{\tau}^{\alpha}((u^n+\chi U)\partial_xu^{n+1}-(f_1^n+F_1)\partial_xf_1^{n+1}-(f_2^n+F_2)\partial_xf_2^{n+1}), \partial_t\partial_\tau^{\alpha}u^{n+1}\rangle}_{K_1}\notag\\
&\quad+\underbrace{ \langle \partial_{\tau}^{\alpha}(\chi\partial_xUu^{n+1}-(\partial_xF_1+\chi^{\prime\prime}UF_1)f_1^{n+1}-\partial_xF_2f_2^{n+1}), \partial_t\partial_\tau^{\alpha}u^{n+1}\rangle}_{K_2}\notag\\
&\quad\underbrace{-\langle \partial_\tau^\alpha\left[f_k^n\partial_y(f_k^n\partial_yu^{n+1}+\chi^\prime Uf_1^{n+1})\right], \partial_t\partial_\tau^{\alpha}u^{n+1}\rangle}_{K_3}=\underbrace{\langle \partial_\tau^\alpha r_1, \partial_t\partial_\tau^{\alpha}u^{n+1}\rangle}_{K_4}.\label{P2.4-1}
\end{align}
By a straightforward computation, we obtain
\begin{align}
|K_1|+|K_2|&\le C(1+\|{\bf\Lambda}^n\|_{H^{m-1}}^2)\|{\bf\Lambda}^{n+1}\|_{H^{m-1}}^2+\frac{1}{4}\|\partial_t\partial_\tau^{\alpha}u^{n+1}\|_{L^2}^2,\label{P2.4-2}\\
|K_4|&\le C\|r_1\|_{m-1}^2+\frac{1}{8}\|\partial_t\partial_\tau^{\alpha}u^{n+1}\|_{L^2}^2\leq \frac{1}{8}\|\partial_t\partial_\tau^{\alpha}u^{n+1}\|_{L^2}^2+C .\label{P2.4-3}
\end{align}

For the diffusion term $K_3$, we have
\begin{align*}
K_3&=-\langle f_k^n\partial_\tau^{\alpha}\partial_y(f_k^{n}\partial_yu^{n+1}), \partial_t\partial_\tau^{\alpha}u^{n+1}\rangle\\
&\quad-\langle \mathcal{C}_k\partial_y(f_k^{n}\partial_yu^{n+1}), \partial_t\partial_\tau^{\alpha}u^{n+1}\rangle\\
&\quad-\langle \partial_\tau^{\alpha}\left[f_k^n\partial_y(\chi^\prime Uf_1^{n+1})\right], \partial_t\partial_\tau^\alpha u^{n+1}\rangle\\
&=:K_{31}+K_{32}+K_{33}.
\end{align*}
Integrating by parts and using the boundary condition $\partial_t\partial_\tau^{\alpha}u^{n+1}|_{y=0}=0$, we deduce that
\begin{align*}
K_{31}&=\langle f_k^n\partial_\tau^{\alpha}(f_k^{n}\partial_yu^{n+1}), \partial_t\partial_\tau^{\alpha}\partial_yu^{n+1}+\langle \partial_yf_k^n\partial_\tau^{\alpha}(f_k^{n}\partial_yu^{n+1}), \partial_t\partial_\tau^{\alpha}u^{n+1}\rangle\\
&=\frac{1}{2}\frac{d}{dt}\|f_k^n\partial_\tau^{\alpha}\partial_yu^{n+1}\|_{L^2}^2-\frac{1}{2}\int_{\Omega}\partial_t(f_k^n)^2(\partial_\tau^{\alpha}\partial_yu^{n+1})^2\;dxdy\nonumber\\
&\quad+\langle f_k^n\mathcal{C}_k\partial_yu^{n+1}, \partial_t\partial_\tau^{\alpha}\partial_yu^{n+1}\rangle+\langle \partial_yf_k^n\partial_\tau^{\alpha}(f_k^{n}\partial_yu^{n+1}), \partial_t\partial_\tau^{\alpha}u^{n+1}\rangle\nonumber.
\end{align*}
Using the Cauchy-Schwarz inequality and the commutator estimate \eqref{com}, we obtain
\begin{align}
|\langle f_k^n\mathcal{C}_k\partial_yu^{n+1}, \partial_t\partial_\tau^{\alpha}\partial_yu^{n+1}\rangle|&\le \delta\|\partial_t\partial_\tau^\alpha\partial_yu^{n+1}\|_{L^2}^2+ C_{\delta}\int_{\Omega} \left[f_k^n\mathcal{C}_k\partial_yu^{n+1}\right]^2\;dxdy
\notag\\
&\le  \delta\|\partial_t\partial_\tau^\alpha\partial_yu^{n+1}\|_{L^2}^2+C_{\delta}(1+\|(f_1^n, f_2^n)\|_{H^{m-1}}^4)\|\partial_yu^{n+1}\|_{H^{m-1}}^2\nonumber\\
&\leq \delta\|\partial_t\partial_\tau^\alpha\partial_yu^{n+1}\|_{L^2}^2+C_{\delta}(1+{\bf \|\Lambda}^n\|_{H^{m-1}}^4)\|{\bf \Lambda}^{n+1}\|_{H^{m}}^2,\nonumber
\end{align}
where $\delta>0$ can be chosen arbitrarily small. Similarly, we have
\begin{align}
|\langle \partial_yf_k^n\partial_\tau^{\alpha}(f_k^{n}\partial_yu^{n+1}), \partial_t\partial_\tau^{\alpha}u^{n+1}\rangle|&\le \frac{1}{8}\|\partial_t\partial_\tau^{\alpha}u^{n+1}\|_{L^2}+C\int_{\Omega} \left[\partial_yf_k^n\partial_\tau^{\alpha}(f_k^{n}\partial_yu^{n+1})\right]^2\;dxdy
\notag\\
&\le \frac{1}{8}\|\partial_t\partial_\tau^{\alpha}u^{n+1}\|_{L^2}^2+C(1+\|(f_1^n, f_2^n)\|_{H^{m-1}}^4)\|\partial_yu^{n+1}\|_{H^{m-1}}^2\nonumber\\
&\leq \frac{1}{8}\|\partial_t\partial_\tau^{\alpha}u^{n+1}\|_{L^2}^2+C(1+{\bf \|\Lambda}^n\|_{H^{m-1}}^4)\|{\bf \Lambda}^{n+1}\|_{H^{m}}^2.\nonumber
\end{align}
The remaining two terms of $K_3$ can be estimated in the following way:
\begin{align}
|K_{32}|+|K_{33}|&\le \frac{1}{4}\|\partial_t\partial_\tau^{\alpha}u^{n+1}\|_{L^2}^2+ C\int_{\Omega} \left[(\mathcal{C}_k\partial_y(f_k^{n}\partial_yu^{n+1})\right]^2+\left[\partial_\tau^{\alpha}(f_k^n\partial_y(\chi^\prime Uf_1^{n+1}))\right]^2\;dxdy\notag\\
&\le \frac{1}{4}\|\partial_t\partial_\tau^{\alpha}u^{n+1}\|_{L^2}^2+ C(1+\|(f_1^n, f_2^n)\|_{H^{m-1}}^4)\|\partial_yu^{n+1}\|_{H^{m-1}}^2\notag\\
&\quad+C(1+\|(f_1^n, f_2^n)\|_{H^{m-1}}^2)\|f_1^{n+1}\|_{H^{m}}^2\nonumber\\
&\leq \frac{1}{4}\|\partial_t\partial_\tau^{\alpha}u^{n+1}\|_{L^2}^2+C(1+{\bf \|\Lambda}^n\|_{H^{m-1}}^4)\|{\bf \Lambda}^{n+1}\|_{H^{m}}^2.\nonumber
\end{align}
Combining all the above inequalities, we arrive at
\begin{align}
K_3\geq &\frac{1}{2}\frac{d}{dt}\|f_k^n\partial_\tau^{\alpha}\partial_yu^{n+1}\|_{L^2}^2-\frac{3}{8}\|\partial_t\partial_\tau^{\alpha}u^{n+1}\|_{L^2}^2\nonumber\\
&-\delta\|\partial_t\partial_\tau^\alpha\partial_yu^{n+1}\|_{L^2}^2-C_\delta(1+{\bf \|\Lambda}^n\|_{H^{m-1}}^4)\|{\bf \Lambda}^{n+1}\|_{H^{m}}^2.\label{P2.4-4}
\end{align}
Substituting \eqref{P2.4-2}, \eqref{P2.4-3} and \eqref{P2.4-4} into \eqref{P2.4-1} yields \eqref{ad3}.
\end{proof}

\subsection{Estimates of higher-order normal derivatives of velocity}
\begin{lemma}\label{P2.5}
For any solutions to \eqref{2.16}, it holds that
\begin{align}
\|\partial_y^2u^{n+1}\|_{H^{m-1}}\le& C(1+\|{\bf\Lambda}^{n}\|_{H^{m}}^2)\|{\bf\Lambda}^{n+1}\|_{H^{m}}+C,\label{ad4}\\
\|\partial_y^2u^{n+1}\|_{H^{m-2}}\le& C+C\int_0^t(1+\|{\bf\Lambda}^{n}(s)\|_{H^{m}}^2)\|{\bf\Lambda}^{n+1}(s)\|_{H^{m}}ds.\label{ad4-1}
\end{align}
\end{lemma}

\begin{proof}
We rewrite the first equation of \eqref{2.16} as follows.
\begin{align}\label{2.36}
(f_k^n)^2\partial_y^2u^{n+1}&=\partial_tu^{n+1}+(u^n+\chi U)\partial_xu^{n+1}-(f_1^n+F_1)\partial_xf_1^{n+1}-(f_2^n+F_2)\partial_xf_2^{n+1}\notag\\
&\quad+\chi\partial_xU u^{n+1}-(\partial_xF_1+\chi^{\prime\prime}UF_1)f_1^{n+1}-\partial_xF_2f_2^{n+1}\notag\\
&\quad-f_k^n\partial_yf_1^n\partial_yu^{n+1}-f_k^n\partial_y(\chi^\prime Uf_1^{n+1})-r_1.
\end{align}
Applying the full derivatives $D^\alpha$ on the both side of \eqref{2.36} with $|\alpha|\le m-1$,  we obtain
\begin{align}
(f_k^n)^2D^{\alpha}\partial_y^2u^{n+1}&=-[D^{\alpha}, (f_k^n)^2]\partial_y^2u^{n+1}+D^{\alpha}\partial_tu^{n+1}-D^{\alpha}r_1+D^\alpha \mathcal{I}_{re},\label{P2.5-1}
\end{align}
where
$$
\begin{aligned}
\mathcal{I}_{re}=&(u^n+\chi U)\partial_xu^{n+1}-(f_1^n+F_1)\partial_xf_1^{n+1}-(f_2^n+F_2)\partial_xf_2^{n+1}+\chi\partial_xU u^{n+1}\\
&-(\partial_xF_1+\chi^{\prime\prime}UF_1)f_1^{n+1}-\partial_xF_2f_2^{n+1}-f_k^n\partial_yf_1^n\partial_yu^{n+1}-f_k^n\partial_y(\chi^\prime Uf_1^{n+1}).
\end{aligned}
$$

For the first term involving the commutator, we have
\begin{align}
\|[D^{\alpha}, (f_k^n)^2]\partial_y^2u^{n+1}\|_{L^2}&\leq C \sum\limits_{|\beta|+|\gamma|=|\alpha| \atop |\beta|\ge 1}\|D^{\beta}(f_k^n)^2D^{\gamma}\partial_y^2u^{n+1}\|_{L^2}\nonumber\\
&\leq C(1+\|(f_1^n, f_2^n)\|_{H^{m-1}}^2)\|\partial_y^2u^{n+1}\|_{H^{m-2}}\nonumber\\
&\leq C(1+\|{\bf\Lambda}^{n}\|_{H^{m}}^2)\|{\bf \Lambda}^{n+1}\|_{H^{m}}.\label{P2.5-2}
\end{align}
For the last term in \eqref{P2.5-1}, we have
\begin{align}
\|D^\alpha \mathcal{I}_{re}\|_{L^2}&\le C(1+\|{\bf \Lambda}^n\|_{H^{m-1}})(\|\partial_x{\bf \Lambda}^{n+1}\|_{H^{m-1}}+\|{\bf \Lambda}^{n+1}\|_{H^{m-1}}+\|\partial_yf_1^{n+1}\|_{H^{m-1}})\notag\\
&\quad+C(1+\|(f_1^n, f_2^n)\|_{H^{m-1}})\|\partial_yf_1^n\|_{H^{m-1}}\|\partial_yu^{n+1}\|_{H^{m-1}}\nonumber\\
&\leq C(1+\|{\bf\Lambda}^{n}\|_{H^{m}}^2)\|{\bf \Lambda}^{n+1}\|_{H^{m}}.\label{P2.5-3}
\end{align}
Substituting \eqref{P2.5-2} and \eqref{P2.5-3} into \eqref{P2.5-1} and using the lower bound \eqref{2.15} of $f_k^n$, we deduce that
\begin{align}
\|\partial_y^2u^{n+1}\|_{H^{m-1}}\le C(1+\|{\bf\Lambda}^{n}\|_{H^{m}}^2)\|{\bf \Lambda}^{n+1}\|_{H^{m}}+C,\nonumber
\end{align}
which is the desired estimate \eqref{ad4}. For \eqref{ad4-1}, note that
$$
\begin{aligned}
\|\partial_y^2u^{n+1}(t)\|_{H^{m-2}}&\leq \|\partial_y^2u^{n+1}(0)\|_{H^{m-2}}+\int_0^t\|\partial_t\partial_y^2u^{n+1}(s)\|_{H^{m-2}}ds\nonumber\\
&\leq \|\partial_y^2u^{n+1}(0)\|_{H^{m-2}}+\int_0^t\|\partial_y^2u^{n+1}(s)\|_{H^{m-1}}ds.
\end{aligned}
$$
Substituting \eqref{ad4} into above equation yields \eqref{ad4-1}. The proof of Lemma \ref{P2.5} is complete.
\end{proof}

\subsection{Estimates of normal derivatives of $(f_1^{n+1},f_2^{n+1})$}
Finally, we estimate the normal derivatives of the elastic part $(f_1^{n+1},f_2^{n+1})$. For this, we need to use the symmetric-hyperbolic structure of the elastic equations.
\begin{lemma}\label{P2.6}
Let $\delta\in (0,1)$, which can be chosen arbitrarily small. There exists a positive constant $C_\delta$, such that
\begin{align}\label{2.60}
\frac{d}{dt}\|(\partial_yf_1^{n+1},\partial_yf^{n+1}_2)\|_{H^{m-1}}^2&\leq \delta\|\partial_yu^{n+1}\|_{H^m}^2+C_{\delta}(1+\|{\bf \Lambda}^n\|_{H^m}^2)\|{\bf \Lambda}^{n+1}\|_{H^m}^2+C_\delta.
\end{align}
\end{lemma}

\begin{proof}
Applying the operator $D^{\alpha}\partial_y$ to the second equation of \eqref{2.16} with $|\alpha|\le m-1$, and taking inner product of the resulting equation with $D^{\alpha}\partial_yf_1^{n+1}$, we obtain
\begin{align}
&\frac{1}{2}\frac{d}{dt}\|D^{\alpha}\partial_yf_1^{n+1}\|_{L^2}^2+
\underbrace{\langle D^{\alpha}\partial_y\left[(u^n+\chi U)\partial_xf_1^{n+1}\right], D^{\alpha}\partial_yf_1^{n+1}}_{L_1}\rangle\nonumber\\
&\quad\underbrace{-\langle D^{\alpha}\partial_y\left[(f_1^n+F_1)\partial_xu^{n+1}\right], D^{\alpha}\partial_yf_1^{n+1}\rangle}_{L_2}
+\underbrace{\langle D^{\alpha}\partial_y(\partial_xF_1u^{n+1}), D^{\alpha}\partial_yf_1^{n+1}}_{L_3}\rangle\nonumber\\
&\quad\underbrace{-\langle D^{\alpha}\partial_y(\chi\partial_xUf_1^{n+1}), D^{\alpha}\partial_yf_1^{n+1}\rangle}_{L
_4}=\underbrace{\langle D^{\alpha}\partial_yr_2, D^{\alpha}\partial_yf_1^{n+1}\rangle}_{L_5}.\label{P2.6-1}
\end{align}
We split $L_1$ into the following two parts,
\begin{align*}
L_1=\langle (u^n+\chi U)\partial_xD^{\alpha}\partial_yf_1^{n+1}, D^{\alpha}\partial_yf_1^{n+1}\rangle+\langle \mathcal{C}_4\partial_xf_1^{n+1}, D^{\alpha}\partial_yf_1^{n+1}\rangle=: L_{11}+L_{12},
\end{align*}
where $\mathcal{C}_4$ is the commutator defined by
\begin{align*}
\mathcal{C}_4=[D^{\alpha}\partial_y, u^n+\chi U].
\end{align*}
For the first term $L_{11}$, we integrate by part to obtain
\begin{align}
|L_{11}|&=\frac{1}{2}\left|-\langle \partial_x(u^n+\chi U)D^{\alpha}\partial_yf_1^{n+1}, D^{\alpha}\partial_yf_1^{n+1}\rangle\right|\notag\\
&\le C(1+\|\partial_xu^n\|_{L^\infty})\|D^{\alpha}\partial_yf_1^{n+1}\|_{L^2}\nonumber\\
&\leq C(1+\|{\bf \Lambda}^n\|_{H^m})\|{\bf \Lambda}^{n+1}\|_{H^m}^2.\label{P2.6-3}
\end{align}
For the commutator term, note that
\begin{align*}
\|\mathcal{C}_4\partial_xf_1^{n+1}\|_{L^2}&\lesssim \sum\limits_{|\beta|+|\gamma|=|\alpha| \atop |\beta|\ge 0}\|D^{\beta}\partial_y(u^n+\chi U)D^{\gamma}\partial_xf_1^{n+1}\|_{L^2}\\
&\quad +\sum\limits_{|\beta|+|\gamma|=|\alpha| \atop |\beta|\ge 1}\|D^{\beta}(u^n+\chi U)D^{\gamma}\partial_y\partial_xf_1^{n+1}\|_{L^2}\nonumber\\
&\leq C(1+\|u^n\|_{H^{m}})\|\partial_xf_1^{n+1}\|_{H^{m-1}}.\nonumber
\end{align*}
Thus, we obtain
\begin{align}\label{P2.6-2}
|L_{12}|\le & C(1+\|u^n\|_{H^{m}})\|\partial_xf_1^{n+1}\|_{H^{m-1}}\|D^{\alpha}\partial_yf_1^{n+1}\|_{L^2}\nonumber\\
\leq & C(1+\|{\bf \Lambda}^n\|_{H^m})\|{\bf \Lambda}^{n+1}\|_{H^m}^2.
\end{align}

Similar ato $L_1$, we write $L_2$ in the following form
\begin{align*}
L_2=\langle (f_1^n+F_1)D^{\alpha}\partial_x\partial_yu^{n+1}, D^{\alpha}\partial_yf_1^{n+1}\rangle+\langle \mathcal{C}_5\partial_xu^{n+1}, D^{\alpha}\partial_yf_1^{n+1}\rangle=: L_{21}+L_{22}.
\end{align*}
Here, $\mathcal{C}_5$ denotes the commutator
\begin{align*}
\mathcal{C}_5=[D^{\alpha}\partial_y, f_1^n+F_1].
\end{align*}
Using the Cauchy-Schwarz inequality, we obtain
\begin{align}
|L_{21}|&\le \delta\|D^{\alpha}\partial_x\partial_yu^{n+1}\|_{L^2}^2+C_{\delta}\int_{\Omega} (f_1^n+F_1)^2(D^{\alpha}\partial_yf_1^{n+1})^2\;dxdy\notag\\
&\le \delta\|D^{\alpha}\partial_x\partial_yu^{n+1}\|_{L^2}^2+ C_{\delta}(1+\|f_1^n\|_{L^\infty}^2)\|D^{\alpha}\partial_yf_1^{n+1}\|_{L^2}^2\nonumber\\
&\leq \delta\|\partial_yu^{n+1}\|_{H^{m}}^2+C_\delta(1+\|{\bf \Lambda}^n\|_{H^m}^2)\|{\bf \Lambda}^{n+1}\|_{H^m}^2,\label{P2.6-4}
\end{align}
where $\delta\in (0,1)$ can be chosen arbitrarily small. Similar as \eqref{P2.6-2}, we control the commutator term $L_{22}$ as follows:
\begin{align}
|L_{22}|&\le C(1+\|f_1^n\|_{H^{m}})\|\partial_xu^{n+1}\|_{H^{m-1}}\|D^{\alpha}\partial_yf_1^{n+1}\|_{L^2}\nonumber\\
&\leq C(1+\|{\bf \Lambda}^n\|_{H^m})\|{\bf \Lambda}^{n+1}\|_{H^m}^2.\label{P2.6-5}
\end{align}

The remaining terms $L_3$ and $L_4$, which respectively involve the zeroth-order and error terms, can be estimated directly:
\begin{align}
|L_3|&\le C\|(\partial_yu^{n+1}, \partial_yf_1^{n+1})\|_{H^{m-1}}\|D^{\alpha}\partial_yf_1^{n+1}\|_{L^2}\leq C\|{\bf \Lambda}^{n+1}\|_{H^m}^2,\label{P2.6-6}\\
|L_4|&\le C\|\partial_yr_2\|_{H^{m-1}}\|D^{\alpha}\partial_yf_1^{n+1}\|_{L^2}\leq C\|{\bf \Lambda}^{n+1}\|_{H^m}^2+C.\label{P2.6-7}
\end{align}
Substituting \eqref{P2.6-3}-\eqref{P2.6-7} into \eqref{P2.6-1} and summing over all $\alpha$ with $|\alpha|\le m-1$, we arrive at
\begin{align}
\frac{d}{dt}\|\partial_yf_1^{n+1}\|_{H^{m-1}}^2\le \delta\|\partial_yu^{n+1}\|_{H^m}^2+C_{\delta}(1+\|{\bf\Lambda}^n\|_{H^m}^2)\|{\bf \Lambda}^{n+1}\|_{H^m}^2+C.\nonumber
\end{align}
The estimate of $f_2^{n+1}$ can be obtained in the same way. The proof of Lemma \ref{P2.6} is complete.
\end{proof}

\subsection{Uniform bounds on $(u^{n},f_1^n,f_2^n)$}
From Lemmas \ref{P2.2}-\ref{P2.6}, we obtain the following differential inequalities:
\begin{align}\label{2.50}
&\frac{d}{dt}\left(\|{\bf \Lambda}^{n+1}\|_{m}^2+\sum\limits_{|\alpha|\le m-1}\|f_k^n\partial_\tau^{\alpha}\partial_yu^{n+1}\|_{L^2}^2+\|(\partial_yf_1^{n+1}, \partial_yf_2^{n+1})\|_{H^{m-1}}^2\right)\notag\\
&\qquad+\|\partial_yu^{n+1}\|_{m}^2+\|\partial_y^2u^{n+1}\|_{H^{m-1}}^2\leq \delta\|\partial_yu^{n+1}\|_{H^m}^2+C_{\delta}(1+\|{\bf \Lambda}^n\|_{H^m}^4)\|{\bf \Lambda}^{n+1}\|_{H^m}^2+C,
\end{align}
where $\delta\in (0,1)$ can be chosen arbitrarily small.  Notice that
\begin{align*}
\|\partial_yu^{n+1}\|_{H^m}\lesssim \|\partial_yu^{n+1}\|_{m}+\|\partial_y^2u^{n+1}\|_{H^{m-1}}.
\end{align*}
Then by taking $\delta$ suitably small, such that the first term on the right-hand side of \eqref{2.50} can be absorbed, we obtain
\begin{align*}
&\frac{d}{dt}\left(\|{\bf\Lambda}^{n+1}\|_{m}^2+\sum\limits_{|\alpha|\le m-1}\|f_k^n\partial_\tau^{\alpha}\partial_yu^{n+1}\|_{L^2}^2+\|(\partial_yf_1^{n+1}, \partial_yf_2^{n+1})\|_{H^{m-1}}^2\right)\notag\\
&\qquad+\|\partial_yu^{n+1}\|_{m}^2+\|\partial_y^2u^{n+1}\|_{H^{m-1}}^2\le C(1+\|{\bf \Lambda}^n\|_{H^m}^4)\|{\bf \Lambda}^{n+1}\|_{H^m}^2+C.
\end{align*}

Next, we integrate it with respect to $t$ and get
\begin{align}\label{2.70}
&\|{\bf \Lambda}^{n+1}\|_{m}^2+\sum\limits_{|\alpha|\le m-1}\|f_k^n\partial_\tau^{\alpha}\partial_yu^{n+1}\|_{L^2}^2+\|(\partial_yf_1^{n+1}, \partial_yf_2^{n+1})\|_{H^{m-1}}^2\notag\\
&\qquad+\int_0^t\|\partial_yu^{n+1}(s)\|_{m}^2+\|\partial_y^2u^{n+1}(s)\|_{H^{m-1}}^2\;ds\notag\\
&\quad\le \|{\bf \Lambda}^{n+1}(0)\|_{m}^2+\sum\limits_{|\alpha|\le m-1}\|f_k^n\partial_\tau^{\alpha}\partial_yu^{n+1}(0)\|_{L^2}^2+\|(\partial_yf_1^{n+1}, \partial_yf_2^{n+1})(0)\|_{H^{m-1}}^2\notag\\
&\qquad+C\int_0^t(1+\|{\bf \Lambda}^n(s)\|_{H^m}^4)\|{\bf \Lambda}^{n+1}(s)\|_{H^m}^2\;ds+C\notag\\
&\quad\le C+C\int_0^t(1+\|{\bf\Lambda}^n(s)\|_{H^m}^4)\|{\bf\Lambda}^{n+1}(s)\|_{H^m}^2\;ds.
\end{align}
Moreover, substituting bounds \eqref{ad4}, \eqref{ad4-1} and \eqref{2.60} on normal derivatives into \eqref{2.70}, we deduce that
\begin{align}
&\|{\bf\Lambda}^{n+1}\|_{H^m}^2+\int_0^t\|\partial_yu^{n+1}(s)\|_{m}^2+\|\partial_y^2u^{n+1}(s)\|_{H^{m-1}}^2\;ds\notag\\
&\quad\lesssim \|{\bf\Lambda}^{n+1}\|_{m}^2+\|\partial_yu^{n+1}\|_{m-1}^2+\|\partial_y^2u^{n+1}\|_{H^{m-2}}+\|(\partial_yf_1^{n+1}, \partial_yf_2^{n+1})\|_{H^{m-1}}^2\notag\\
&\qquad+\int_0^t\|\partial_yu^{n+1}(s)\|_{m}^2+\|\partial_y^2u^{n+1}(s)\|_{H^{m-1}}^2\;ds\notag\\
&\quad\le C+\int_0^t C(1+\|{\bf\Lambda}^{n}(s)\|_{H^{m}}^4)\|{\bf\Lambda}^{n+1}(s)\|_{H^{m}}^2\;ds.\nonumber
\end{align}
Then by the Gronwall inequality, there exists a positive constant $C_1$, such that
\begin{align}\label{2.68}
\|{\bf\Lambda}^{n+1}\|_{H^m}^2+\int_0^t\|\partial_yu^{n+1}(s)\|_{m}^2+\|\partial_y^2u^{n+1}(s)\|_{H^{m-1}}^2\;ds\le C_1\exp\left\{ C_1\int_0^t \|{\bf\Lambda}^{n}(s)\|_{H^{m}}^4\;ds\right\}.
\end{align}

From the Gronwall-type estimate \eqref{2.68}, we obtain the following uniform bounds of approximate solutions $(u^{n},f^n_1,f^n_2)$.
\begin{proposition}\label{prop1}
Let
\begin{align}
C_0=\max\{ C_1, \sup\limits_{t\in [0, T_0]}\|\Lambda^0(t)\|_{H^m}\}\nonumber
\end{align}
where $C_1$ is given by \eqref{2.68} and set $T_*=\min\{T_0, \frac{15}{32C_0^3},\hat{T}\}$, where $\hat{T}$ is given by
\begin{align}
\hat{T}(|k_1|+|k_2|)\left(2C_0^{\frac12}+M\right)=\sigma_0,\nonumber
\end{align}
and $M$ as given in \eqref{2.10}.
Then the approximate solution sequence $\{{\bf\Lambda}^n\}(n\ge 0)$ satisfies
\begin{align}\label{2.76}
\|{\bf\Lambda}^n\|_{H^m}^2\le \frac{C_0}{\sqrt{1-2C_0^3t}},
\end{align}
and
\begin{align}
f_k^{n+1}(t, x, y)\ge \sigma_0,\label{2.76-1}
\end{align}
for all $(t, x, y)\in (0, T_*)\times \mathbb{T}\times \mathbb{R}_+$.
\end{proposition}
\begin{proof}
It is straightforward to check \eqref{2.76} and \eqref{2.76-1} holds for $n=0$. Assume that they holds for some $n>0$, then from \eqref{2.68}, we deduce that
\begin{align}
\|{\bf \Lambda}^{n+1}\|_{H^m}^2&\le C_0\exp\left\{ C_0\int_0^t \|{\bf\Lambda}^{n}(s)\|_{H^{m}}^4\;ds\right\}\notag\\
&\le C_0\exp\left\{ \int_0^t\frac{C_0^3}{1-2C_0^3s}\;ds\right\}=\frac{C_0}{\sqrt{1-2C_0^3t}}.\nonumber
\end{align}
Moreover, for any $t\in [0,T_*]$, we have
\begin{align*}
f_k^{n+1}(t, x, y)&=k_1(f_1^{n+1}(0, x, y)+F_1(0, x))+k_2(f_2^{n+1}(0, x, y)+F_2(0, x))\\
&\quad+k_1\int_0^t \partial_t(f_1^{n+1}(s, x, y)+F_1(s, x))\;ds+k_2\int_0^t \partial_t(f_2^{n+1}(s, x, y)+F_2(s, x))\;ds\\
&\ge 2\sigma_0-|k_1|\int_0^t \|\partial_tf_1^{n+1}(s, \cdot)\|_{L^\infty}+\|\partial_tF_1(s, \cdot)\|_{L^\infty}\;ds\\
&\quad\qquad-|k_2|\int_0^t \|\partial_tf_2^{n+1}(s, \cdot)\|_{L^\infty}+\|\partial_tF_2(s, \cdot)\|_{L^\infty}\;ds\\
&\ge 2\sigma_0-|k_1|\int_0^t \|f_1^{n+1}\|_{H^3}+\|\partial_tF_1\|_{H^1}\;ds-|k_2|\int_0^t \|f_2^{n+1}\|_{H^3}+\|\partial_tF_2\|_{H^1}\;ds\\
&\ge 2\sigma_0-t(|k_1|+|k_2|)\left(\frac{C_0^{\frac12}}{(1-2C_0^3t)^{\frac14}}+M\right)\nonumber\\
&\ge 2\sigma_0-t(|k_1|+|k_2|)\left(2C_0^{\frac12}+M\right)\geq \sigma_0.
\end{align*}
Therefore, \eqref{2.76} and \eqref{2.76-1} follow from the induction argument. The proof of Proposition \ref{prop1} is complete.
\end{proof}

\subsection{Proof of Proposition \ref{T2.1} }
We first prove the approximation sequence $(u^n,f^n_1,f_2^n)$ is a Cauchy sequence in $L^\infty(0,T_*; L^2(\Omega))$. For any $n\ge 1$, we introduce the difference:
\begin{align*}
\tilde{\bf{\Lambda}}^{n}=(\tilde{u}^{n}, \tilde{f}_1^n, \bar{f}_2^n)={\bf\Lambda}^{n+1}-{\bf\Lambda}^n=(u^{n+1}-u^n, f_1^{n+1}-f_1^n, f_2^{n+1}-f_2^n),
\end{align*}
and
\begin{align*}
\tilde{f}_k^n=k_1\tilde{f}_1^n+k_2\tilde{f}_2^n,
\end{align*}
where $(k_1,k_2)$ is defined in Theorem \ref{T1.1}. It is straightforward to check that $(\tilde{u}^{n}, \tilde{f}_1^n, \tilde{f}_2^n)$ satisfies:
\begin{align}
\begin{cases}
\partial_t\tilde{u}^{n}+(u^n+\chi U)\partial_x\tilde{u}^{n}-f_k^n\partial_y(f_k^n\partial_y\tilde{u}^{n}+\chi^\prime U\tilde{f}_1^{n})-(f_1^n+F_1)\partial_x\tilde{f}_1^{n}\\
\qquad-(f_2^n+F_2)\partial_x\tilde{f}_2^{n}+\chi\partial_xU \tilde{u}^{n}-(\partial_xF_1+\chi^{\prime\prime}UF_1)\tilde{f}_1^{n}-\partial_xF_2\tilde{f}_2^{n}-f_k^{n-1}\partial_y(\tilde{f}_k^n\partial_yu^n)=\tilde{r}_1^n,\\
\partial_t\tilde{f}_1^{n}+(u^n+\chi U)\partial_x\tilde{f}_1^{n}-(f_1^n+F_1)\partial_x\tilde{u}^{n}+\partial_xF_1\tilde{u}^{n}-\chi\partial_xU\tilde{f}_1^{n}=\tilde{r}_2^n,\\
\partial_t\tilde{f}_2^{n}+(u^n+\chi U)\partial_x\tilde{f}_2^{n}-(f_2^n+F_2)\partial_x\tilde{u}^{n}+\partial_xF_2\tilde{u}^{n}-\chi\partial_xU\tilde{f}_2^{n}=\tilde{r}_3^n,
\end{cases}\nonumber
\end{align}
where the error terms $\tilde{r}_1^n, \tilde{r}_2^n, \tilde{r}_3^n$ are defined by
\begin{align*}
\tilde{r}_1^n&=-\tilde{u}^{n-1}\partial_xu^n+\tilde{f}_1^{n-1}\partial_xf_1^n+\tilde{f}_2^{n-1}\partial_xf_2^n
+\tilde{f}_k^{n-1}\partial_y(f_k^n\partial_yu^n+\chi^\prime Uf_1^n),\\
\tilde{r}_2^n&=-\tilde{u}^{n-1}\partial_xf_1^{n}+\tilde{f}_1^{n-1}\partial_xu^n,\\
\tilde{r}_3^n&=-\tilde{u}^{n-1}\partial_xf_2^{n}+\tilde{f}_2^{n-1}\partial_xu^n.
\end{align*}

Applying the same approach as in Lemma \ref{P2.2} to $\tilde{{\bf \Lambda}}^n$  , we deduce the following $L^2$-estimate:
\begin{align}
\frac{d}{dt}\|\tilde{\bf \Lambda}^{n}\|_{L^2}^2+\|\partial_y\tilde{u}^{n}\|_{L^2}^2&\le C(1+\|{\bf\Lambda}^n\|_{H^{3}}^2)\|\tilde{\bf \Lambda}^{n}\|_{L^2}^2+C\left|\int_{\Omega} f_k^{n-1}\partial_y(\tilde{f}_k^n\partial_yu^n)\cdot\tilde{u}^{n}\;dxdy\right|\notag\\
&\quad+C\|(\tilde{r}_1^n,\tilde{r}_2^n, \tilde{r}_3^n)\|_{L^2}\|(\tilde{u}^n,\tilde{f}_1^{n}, \tilde{f}_2^n)\|_{L^2}.\label{2.77-1}
\end{align}
Integrating by parts and using the boundary condition $\tilde{u}^n|_{y=0}=0$, we get
\begin{align*}
&\int_{\Omega} f_k^{n-1}\partial_y(\tilde{f}_k^n\partial_yu^n)\cdot\tilde{u}^{n}\;dxdy\\
&\quad=-\int_{\Omega} f_k^{n-1}\tilde{f}_k^n\partial_yu^n\cdot\partial_y\tilde{u}^{n}\;dxdy-\int_{\Omega} \partial_yf_k^{n-1}\tilde{f}_k^n\partial_yu^n\cdot\tilde{u}^{n}\;dxdy\\
&\quad\le \delta\|\partial_y\tilde{u}^n\|_{L^2}^2+C_\delta\|f_k^{n-1}\partial_yu^n\|_{L^\infty}^2\|\tilde{f}_k^n\|_{L^2}^2+\|\partial_yf_k^{n-1}\partial_yu^n\|_{L^\infty}^2\|(\tilde{f}_k^n, \tilde{u}^n)\|_{L^2}^2\\
&\quad\leq \delta\|\partial_y\tilde{u}^n\|_{L^2}^2+C_{\delta}\|{\bf \Lambda}^n\|_{H^m}^2\|{\bf \tilde{\Lambda}}^n\|_{L^2}^2\nonumber\\
&\quad\leq  \delta\|\partial_y\tilde{u}^n\|_{L^2}^2+C_{\delta}\|{\bf \tilde{\Lambda}}^n\|_{L^2}^2,
\end{align*}
where we have used the uniform estimate \eqref{2.76} in the last inequality. For the remaining part, we note that
\begin{align}
\|(\tilde{r}_{1}^n, \tilde{r}_2^n, \tilde{r}_3^n)\|_{L^2}&\le C\|\partial_x{\bf \Lambda}^n\|_{L^\infty}\|\tilde{\bf \Lambda}^{n-1}\|_{L^2}+C(1+\|{\bf\Lambda}^{n}\|_{H^4}^2)\|(\tilde{f}_1^{n-1}, \tilde{f}_2^{n-1})\|_{L^2}.\nonumber\\
&\leq C(1+\|{\bf \Lambda}^n\|_{H^m}^2)\|{\bf \tilde{\Lambda}}^{n-1}\|_{L^2}\leq C\|{\bf \tilde{\Lambda}}^{n-1}\|_{L^2}.\nonumber
\end{align}

Substituting these bounds into \eqref{2.77-1} and taking $\delta\in (0,1)$ suitably small, we obtain that
\begin{align}
\frac{d}{dt}\|\tilde{\bf \Lambda}^{n}\|_{L^2}^2+\|\partial_y\tilde{u}^{n}\|_{L^2}^2\le C_2\|\tilde{\bf \Lambda}^{n-1}\|_{L^2}^2+C_2\|\tilde{\bf \Lambda}^{n}\|_{L^2}^2\nonumber
\end{align}
for some constants $C_2>0$.
By Gronwall inequality, it holds for any $t\in [0, T_*]$ that
\begin{align}
\sup\limits_{0\le s\le t}\|\tilde{\bf \Lambda}^{n}(s)\|_{L^2}^2&\le C_2e^{C_2t}\int_0^t \|\tilde{\bf \Lambda}^{n-1}(s)\|_{L^2}^2\;ds\notag\\
&\le C_2te^{C_2t}\sup\limits_{0\le s\le t}\|\tilde{\bf \Lambda}^{n-1}(s)\|_{L^2}^2.\nonumber
\end{align}

Taking $T_*$ suitably small such that $C_2T_*e^{C_2T_*}\le \frac{1}{2}$, we obtain
\begin{align}
\sup\limits_{0\le s\le t}\|\tilde{\bf \Lambda}^{n}(s)\|_{L^2}^2\le \frac12\sup\limits_{0\le s\le t}\|\tilde{\bf \Lambda}^{n-1}(s)\|_{L^2}^2, \quad n\ge 1.\nonumber
\end{align}
Therefore, $(u^n,f_1^n,f_2^n)$ is a Cauchy sequence in $L^{\infty}(0,T;L^2(\Omega))$.
Note that this sequence also satisfies a uniform $H^m$ bound \eqref{2.76}. Then there exists ${\bf \Lambda}=(u, f_1, f_2)\in H^m$ such that
\begin{align*}
\lim\limits_{n\to +\infty}{\bf\Lambda}^{n}={\bf\Lambda} ~\mbox{in} \,\, L^\infty(0,T; H^k(\Omega)),~~k<m.
\end{align*}
Moreover, it is straightforward to check that the limit $(u, f_1, f_2)$ is a solution to \eqref{2.9} and satisfies the bound
\begin{align}
\sup_{0\leq t\leq T}\|{\bf\Lambda}(t)\|_{H^k}\le 4C_0.\nonumber
\end{align}
The proof of Proposition \ref{T2.1} is complete.

\subsection{Proof of Theorem \ref{T1.1}}
Let $(\hat{u},\hat{f}_1,\hat{f}_2)$ be the solution constructed in Proposition \ref{T2.1}, and $\phi$ be the stream function defined in \eqref{2.1}. Then we introduce
\begin{align}
{u}(t, x, y)&=\hat{u}\left(t, x, \phi(t,x,y)\right)+U(t, x)\chi(\phi(t,x,y)),\nonumber\\
{f}_1(t,x,y)&=\hat{f}_1\left(t, x, \phi(t,x,y)\right)+F_1(t, x)\chi(\phi(t,x,y)),\nonumber\\
f_2(t, x, y)&=\hat{f}_2\left(t, x, \phi(t,x,y)\right)+F_2(t, x)\chi(\phi(t,x,y))\nonumber,
\end{align}
and $v(t,x,y), f_3(t,x,y),f_4(t,x,y)$ are determined by divergence free conditions in \eqref{1.5}. It is straightforward to check that $({u}, {v}, {f}_1, {f}_2, {f}_3,{f}_4)$ is a solution to \eqref{1.1}. Moreover, from \eqref{b} we deduce that the solution $(u,f_1,f_2)$ satisfies \eqref{b1}. The proof of Theorem \ref{1.1} is complete.

\bigskip

\section{Ill-Posedness}\label{S3}

\subsection{Instability mechanism}\label{S3.1} To illustrate the instability mechanism, we take the time-independent shear flow $(U_s,0,F_{1,s},F_{2,s},0,0)(y)$ as a background solution. The linearization of \eqref{1.1} around this shear flow reads
\begin{align}\label{i1}
\begin{cases}
\partial_tu+U_s\partial_xu+v\partial_yU_s-\partial_y^2u-F_{1, s}\partial_xf_1-f_3F_{1, s}^{\prime}-F_{2, s}\partial_xf_2-f_4F_{2, s}^{\prime}=0,\\
\partial_tf_1+U_s\partial_xf_1+vF_{1, s}^{\prime}-F_{1, s}\partial_xu-f_3\partial_yU_s=0,\\
\partial_tf_2+U_s\partial_xf_2+vF_{2, s}^{\prime}-F_{2, s}\partial_xu-f_4\partial_yU_s=0,\\
\partial_xu+\partial_yv=0,\quad \partial_xf_1+\partial_yf_3=0,\quad \partial_xf_2+\partial_yf_4=0,\\
(u, v, f_1, f_2,f_3,f_4)|_{y=0}={\bf{0}},\quad \lim\limits_{y\to\infty}(u, f_1, f_2)(t, x, y)={\bf{0}}.
\end{cases}
\end{align}
As in \cite{GD10}, we look for the following wave-type solution to \eqref{i1}:
\begin{equation}
\begin{aligned}
(u,v)&=e^{i\eps^{-1}(x+\omega(\eps)t)}(\partial_y\Phi,-i\eps^{-1}\Phi),\\
(f_1,f_3)&=e^{i\eps^{-1}(x+\omega(\eps)t)}(\partial_y\phi,-i\eps^{-1}\phi),\\
(f_2,f_4)&=e^{i\eps^{-1}(x+\omega(\eps)t)}(\partial_y\psi,-i\eps^{-1}\psi).\label{i2}
\end{aligned}
\end{equation}
Here $\Phi, \phi,\psi$ are respectively the stream functions of $(u,v)$, $(f_1,f_3)$ and $(f_2,f_4)$ due to the divergence-free conditions in \eqref{i1}. We consider the regime where $\eps\ll1$, which corresponds to high frequencies.

Substituting \eqref{i2} into \eqref{i1}, we derive the following equations for $(\Phi,\phi,\psi)$:
\begin{equation}
\left\{
\begin{aligned}
&\left(U_s+\omega(\eps)\right)\Phi'-U_s'\Phi+i\eps\Phi'''-F_{1,s}\phi'+F_{1,s}'\phi-F_{2,s}\psi'+F_{2,s}'\psi=0,~y>0,\\
&\left(U_s+\omega(\eps)\right)\phi-F_{1,s}\Phi=0,\\
&\left(U_s+\omega(\eps)\right)\psi-F_{2,s}\Phi=0,\label{i3}
\end{aligned}
\right.
\end{equation}
supplemented with the following boundary conditions:
\begin{align}
\Phi|_{y=0}=\Phi'|_{y=0}=\phi|_{y=0}=\phi'|_{y=0}=\psi|_{y=0}=\psi'|_{y=0}=0.\label{i4}
\end{align}

When $F_{1,s}=F_{2,s}\equiv 0$, \eqref{i3} with \eqref{i4} can be reduced to the following eigenvalue problem for the Prandtl equation:
\begin{equation}
\left\{
\begin{aligned}
&\left(U_s+\omega(\eps)\right)\Phi'-U_s'\Phi+i\eps\Phi'''=0,~y>0,\\
&\Phi|_{y=0}=\Phi'|_{y=0}=0.
\label{i5}
\end{aligned}
\right.
\end{equation}

In their seminal work \cite{GD10}, G\'erard-Varet and Dormy constructed an approximate growing mode of \eqref{i5}  under the assumption that $U_s$ has a unique non-degenerate critical point at $y=a$ satisfying $U_s'(a)=0,U_s''(a)\neq 0)$. Their construction is as follows:
\begin{align}
\omega(\eps)&=-U_s(a)+\eps^{\frac12}\tau,\nonumber\\
\Phi_{app}(y)&=\underbrace{H(y-a)\left(U_s(y)-U_s(a)+\eps^{\frac12}\tau\right)}_{\Phi_{in}}+\underbrace{\eps^{\frac12}V\left(\frac{y-a}{\eps^{\frac14}}\right)}_{\Phi_{sl}}.\nonumber
\end{align}
Here $\tau\in\mathbb{C}$ and $H(\cdot)$ is the classical Heaviside function. Note that $\Phi_{in}$ solves the equation \eqref{i5} up to $O(\eps)$. Both $\Phi_{in}$ and $\Phi''_{in}$ have jumps at $y=a$:
\begin{align}
[\Phi_{in}]_{y=a}=\eps^{\frac12}\tau,~~[\Phi_{in}'']_{y=a}=U_s''(a).\nonumber
\end{align}

The shear layer part $\Phi_{sl}$ is defined in the vicinity of $y=a$ and is used to eliminate jumps generated by $\Phi_{in}$.  It satisfies the following profile equation:
\begin{equation}\label{i8}
\left\{
\begin{aligned}
&\left(\tau+U''(a)\frac{z^2}{2}\right)\frac{dV}{dz}-U_{s}''(a)zV+i\frac{d^3V}{dz^3}=0,~z\neq 0,\\
&[V]_{z=0}=-\tau,~~[\frac{dV}{dz}]_{z=0}=0,~~[\frac{d^2V}{dz^2}]|_{z=0}=-U_s''(a),\\
&\lim_{z\rightarrow \pm \infty}V(z)=0.
\end{aligned}
\right.
\end{equation}

In \cite{GD10}, it is shown that there exists $\tau\in \mathbb{C}$ with $\text{Im}\tau<0$, such that \eqref{i8} admits a non-trivial solution $V$ decaying exponentially at $z=\pm\infty.$ This leads to an instability of the Prandtl equation.

Now we take $\Phi_{app}$ as the leading-order instability profile for the viscoelastic model \eqref{i1}. To this end, we define the approximate elastic stream functions $\phi_{app}$ and $\psi_{app}$ by
\begin{align}
\phi_{app}(y)=\frac{F_{1,s}(y)\Phi_{app}(y)}{U_s(y)+\omega(\eps)},~~\psi_{app}(y)=\frac{F_{2,s}(y)\Phi_{app}(y)}{U_s(y)+\omega(\eps)}.\label{i9}
\end{align}
These expressions follow directly from $\eqref{i3}_2$ and $\eqref{i3}_3$ respectively. It is straightforward to check that $(\Phi_{app},\phi_{app},\psi_{app})\in W^{3,\infty}_\alpha(\mathbb{R}_+).$ Then the approximate growing mode
$$(u_{app},v_{app},f_{1,app},f_{2,app},f_{3,app},f_{4,app})$$ are defined by substituting $(\Phi_{app},\phi_{app},\psi_{app})$ into \eqref{i2}.

Inserting this approximate solution into \eqref{i1}, we deduce that
\begin{align}\label{i10}
\begin{cases}
\partial_tu_{app}+U_s\partial_xu_{app}+v_{app}\partial_yU_s-\partial_y^2u_{app}\\
\qquad-F_{1, s}\partial_xf_{1,app}-f_{3,app}F_{1, s}^{\prime}-F_{2, s}\partial_xf_{2,app}-f_{4,app}F_{2, s}^{\prime}
=e^{i\eps^{-1}(x+\omega(\eps)t)}R(y),\\
\partial_tf_{1,app}+U_s\partial_xf_{1,app}+v_{app}F_{1, s}^{\prime}-F_{1, s}\partial_xu_{app}-f_{3,app}\partial_yU_s=0,\\
\partial_tf_{2,app}+U_s\partial_xf_{2,app}+v_{app}F_{2, s}^{\prime}-F_{2, s}\partial_xu_{app}-f_{4,app}\partial_yU_s=0,\\
\partial_xu_{app}+\partial_yv_{app}=0,\quad \partial_xf_{1,app}+\partial_yf_{3,app}=0,\quad \partial_xf_{2,app}+\partial_yf_{4,app}=0,\\
(u_{app}, v_{app}, f_{1,app}, f_{2,app},f_{3,app},f_{4,app})|_{y=0}={\bf{0}},\quad \lim\limits_{y\to\infty}(u_{app}, f_{1,app}, f_{2,app})(t, x, y)={\bf{0}}.
\end{cases}
\end{align}
Here, the error term $R(y)$ is given by
\begin{align}
R(y)=&-\eps^{-1}\left(U_s(y)-\frac{U_s''(a)}{2}(y-a)^2\right)\partial_y\Phi_{sl}+\eps^{-1}\left(U_s'(y)-U_s''(a)(y-a)\right)\Phi_{sl}\nonumber\\
&+-i\partial_y^3\Phi_{in}+\eps^{-1}\left(F_{1,s}^2+F_{2,s}^2\right)\left(\frac{\partial_y\Phi_{sl}}{U_s+\omega}-\frac{U_s'\Phi_{sl}}{(U_s+\omega)^2}\right),\nonumber
\end{align}
where we have used the fact $(U_s+\omega)\Phi_{in}-U_s'\Phi_{in}=0$ in the derivation of $R(y)$.  Our degeneracy condition $F_{1,s}(a)=F_{1,s}'(a)=F_{2,s}(a)=F_{2,s}'(a)=0$ implies
$$F_{i,s}(y)\sim O(1)|y-a|^2,~~ F_{1,s}^2+F_{2,s}^2\sim O(1)|y-a|^4.
$$
Thus, by using $|U_s+\omega|\gtrsim\eps^{\frac12}$, we deduce that
$$
|R(y)|\leq C\left(\eps^{-1}|y-a|^3+\eps^{-\frac32}|y-a|^4|\right)|\partial_y\Phi_{sl}|+C\left(\eps^{-1}|y-a|^2+\eps^{-2}|y-a|^5|\right)|\Phi_{sl}|+C
$$
in the vicinity of $y=a$ and $R(y)$ exponentially decay away from $y=a$. Therefore, the residual satisfies
\begin{align}
\|R\|_{L^\infty_\alpha}\leq C\eps^{-\frac14}.\label{i11}
\end{align}

\subsection{Proof of Theorem \ref{T1.2}} Now we prove Theorem \ref{T1.2} by a contradiction argument. Assume that \eqref{1.8} is not true, namely, there exists $\delta_0>0$, such that for every $\lambda>0$, it holds
\begin{align}\label{3.16}
\|\mathcal{S}_n(t)\|_{L^{2,\infty}_\alpha(\Omega)\rightarrow L^{2,\infty}_\alpha(\Omega)}\leq Cn^{\mu}e^{\lambda n^{\frac12}t},~~\mu\in [0,\frac14),
\end{align}
for any $t\in [0,\delta_0].$ Then we consider the solution
\begin{align}
(u,f_1,f_2)(t)=\mathcal{S}(t)\left(u_{app}(0),f_{1,app}(0),f_{2,app}(0)\right),\nonumber
\end{align}
where $(u_{app}(t),f_{1,app}(t),f_{2,app}(t))$ is the approximate growing mode constructed in Section \ref{S3.1}.

We now define $(\tilde{u},\tilde{f}_1,\tilde{f}_2)(t,y)=e^{-i\eps^{-1}x}(u(t),f_{1}(t),f_{2}(t))$. Without loss of generality, we assume that $\eps^{-1}$ is an integer. Then it holds that
\begin{align}
(\tilde{u},\tilde{f}_1,\tilde{f}_2)(t,y)=e^{-i\eps^{-1}x}\mathcal{S}_{\eps^{-1}}(t)\left(u_{app}(0),f_{1,app}(0),f_{2,app}(0)\right)\nonumber.
\end{align}

On the one hand, we obtain from \eqref{3.16} that
\begin{align}
\|(\tilde{u},\tilde{f}_1,\tilde{f}_2)(t)\|_{L^\infty_\alpha}\leq C_\lambda\eps^{-\mu}e^{\lambda \eps^{-\frac12}t},~~\mu\in [0,\frac14),\label{3.17}
\end{align}
for any $\lambda>0$, and any $t\in [0,\delta_0]$. On the other hand, since $(u_{app}(t),f_{1,app}(t),f_{2,app}(t))$ satisfies \eqref{i10}, we use Duhamel principle to rewrite the real solution $(\tilde{u},\tilde{f}_1,\tilde{f}_2)(t)$ as
\begin{align}
(\tilde{u},\tilde{f}_1,\tilde{f}_2)(t)=\left(\tilde{u}_{app}(t),\tilde{f}_{1,app}(t),\tilde{f}_{2,app}(t)\right)+e^{-i\eps^{-1}x}\int_0^t\mathcal{S}_{\eps^{-1}}(t-s)\left(e^{i\eps^{-1}(x+\omega(\eps)s)}R(y),0,0\right),\label{3.18}
\end{align}
where $(\tilde{u}_{app}(t),\tilde{f}_{1,app}(t),\tilde{f}_{2,app}(t))=e^{i\eps^{-1}x}({u}_{app}(t),{f}_{1,app}(t),{f}_{2,app}(t))$.

Denote $\lambda_0=|\text{Im}\omega(\eps)|$. Then from construction, we have
\begin{align}
\|\left(\tilde{u}_{app}(t),\tilde{f}_{1,app}(t),\tilde{f}_{2,app}(t)\right)\|_{L^\infty_\alpha}=e^{\lambda_0\eps^{-\frac12}t}.\label{3.19}
\end{align}
From \eqref{3.16} and the estimate \eqref{i11} on the residual $R$, we obtain
\begin{align}
\|e^{-i\eps^{-1}x}\int_0^t\mathcal{S}_{\eps^{-1}}(t-s)\left(e^{i\eps^{-1}(x+\omega(\eps)s)}R(y),0,0\right)\|_{L^\infty_\alpha}&\leq C\eps^{-\mu-\frac14}\int_0^te^{\lambda\eps^{-\frac12}(t-s)}e^{\lambda_0\eps^{-\frac12}s}ds\nonumber\\
&\leq \frac{C\eps^{\frac14-\mu}}{\lambda_0-\lambda}e^{\lambda_0\eps^{-\frac12}t},\label{3.20}
\end{align}
where we have chosen $\lambda<\lambda_0$ in the last inequality. Substituting \eqref{3.19} and \eqref{3.20} into \eqref{3.18}, we deduce that
\begin{align}\label{3.21}
\|(\tilde{u},\tilde{f}_1,\tilde{f}_2)(t)\|_{L^\infty_\alpha}\geq (1-C\eps^{\frac14-\mu})e^{\lambda_0\eps^{-\frac12}t}\geq \frac{1}{2}e^{\lambda_0\eps^{-\frac12}t},
\end{align}
for $\eps$ sufficiently small.

Now we define $T^\eps=\eps^{\frac12}\log\left|\frac{\eps^\mu}{4C_\lambda(\lambda_0-\lambda)}\right|$. Then from \eqref{3.21}, we obtain that
$$\|(\tilde{u},\tilde{f}_1,\tilde{f}_2)(t)\|_{L^\infty_\alpha}\geq 2C_\lambda\eps^{-\mu}e^{\lambda\eps^{-\frac12}t},
$$
which contradicts to \eqref{3.17}. The proof of Theorem \ref{T1.2} is complete.

\bigskip

\appendix
\section{Derivation of the Viscoelastic Boundary Layer Model}\label{A}

In the appendix, we formally derive the viscoelastic boundary layer model \eqref{1.1}. We start with the following two-dimensional incompressible viscoelastic equations in the half plane $\{(x,y)\mid x\in \mathbb{T}, y>0\}$:
\begin{align}
\label{A.1}
\begin{cases}
\partial_t{\bf{u}}^{\varepsilon}+({\bf{u}}^{\varepsilon}\cdot\nabla){\bf{u}}^{\varepsilon}+\nabla p^{\varepsilon}-\varepsilon\Delta{\bf{u}}^{\varepsilon}=\div({\bf{F}}^{\varepsilon}({\bf{F}}^{\varepsilon})^T),\\
\partial_t{\bf{F}}^{\varepsilon}+({\bf{u}}^{\varepsilon}\cdot\nabla){\bf{F}}^{\varepsilon}-\nabla {\bf{u}}^{\varepsilon}\cdot {\bf{F}}^{\varepsilon}=0,\\
\partial_xu^{\varepsilon}+\partial_yv^{\varepsilon}=0,\quad \partial_xf_1^{\varepsilon}+\partial_yf_3^{\varepsilon}=0,\quad \partial_xf_2^{\varepsilon}+\partial_yf_4^{\varepsilon}=0,
\end{cases}
\end{align}
where ${\bf{u}}^{\varepsilon}=(u^{\varepsilon}, v^{\varepsilon})$ represents the velocity field and ${\bf{F}}^{\varepsilon}=\left( \begin{array}{cc} f_1^{\varepsilon}& f_2^{\varepsilon}\\ f_3^{\varepsilon}&  f_4^{\varepsilon} \end{array}\right)$ represents the deformation tensor matrix. We supplement \eqref{A.1} with the following initial data:
\begin{align}
\label{A.2}
{\bf{u}}^{\varepsilon}(0, x, y)={\bf{u}}_0(x, y),\quad {\bf{F}}^{\varepsilon}(0, x, y)={\bf{F}}_0(x, y),
\end{align}
and the following no-slip boundary condition for the velocity field:
\begin{align}
\label{A.3}
{\bf{u}}^{\varepsilon}|_{y=0}={\bf{0}}.
\end{align}
We emphasize that there is no need to impose any boundary condition for ${\bf{F}}^{\varepsilon}$ since its boundary traces can be determined by the initial data. To see this fact, we
evaluate $\eqref{A.1}_2$ at $y=0$ and use the no-slip boundary condition \eqref{A.3} to obtain
\begin{align}
&\partial_t{f}_1^\vep(t,x,0)-{f}^\vep_3(t,x,0)\partial_yu^\vep(t,x,0)=0,\label{A.3-1}\\
&\partial_t{f}_2^\vep(t,x,0)-{f}^\vep_4(t,x,0)\partial_yu^\vep(t,x,0)=0,\label{A.3-2}\\
&\partial_t{f}_3^\vep(t,x,0)=\partial_t{f}_4^\vep(t,x,0)=0.\label{A.3-3}
\end{align}
Thus, from \eqref{A.3-3} we obtain
\begin{align}
&{f}_3^\varepsilon(t,x,0)\equiv {f}_3^{\varepsilon}(0, x, 0)=f_{3, 0}(x, 0),\label{A.4}\\
&{f}_4^\varepsilon(t,x,0)\equiv {f}_4^{\varepsilon}(0, x,0)=f_{4, 0}(x, 0).\label{A.5}
\end{align}
The boundary values $f_{1}^{\vep}(t,x,0)$ and $f_2^{\vep}(t,x,0)$ can be determined by \eqref{A.3-1} and \eqref{A.3-2}.

When the viscosity $\varepsilon$ tends to zero, the system \eqref{A.1} is expected to converge to the following ideal elastodynamic equations
\begin{align}
\label{A.6}
\begin{cases}
\partial_t{\bf{u}}^{0}+({\bf{u}}^{0}\cdot\nabla){\bf{u}}^{0}+\nabla p^{0}=\div({\bf{F}}^{0}({\bf{F}}^{0})^T),\\
\partial_t{\bf{F}}^{0}+({\bf{u}}^{0}\cdot\nabla){\bf{F}}^{0}-\nabla {\bf{u}}^{0}\cdot {\bf{F}}^{0}=0,\\
\partial_xu^{0}+\partial_yv^{0}=0,\quad \partial_xf_1^{0}+\partial_yf_3^{0}=0,\quad \partial_xf_2^{0}+\partial_yf_4^{0}=0,
\end{cases}
\end{align}
with the following non-penetration condition:
\begin{align}
\label{A.7}
v^0|_{y=0}=0.
\end{align}
We supplement \eqref{A.6} with the same initial data as \eqref{A.2}:
\begin{align}
\label{A.8}
{\bf{u}}^{0}(0, x, y)={\bf{u}}_0(x, y),\quad {\bf{F}}^{0}(0, x, y)={\bf{F}}_0(x, y).
\end{align}

It is known from \cite{WX23} that, if a boundary layer forms in the inviscid limit,  the boundary deformation matrix  ${\bf{F}}^\vep|_{y=0}$ must degenerate at the initial time.  Therefore, for simplicity and without loss of generality, we assume that
\begin{align*}
f_{3, 0}(x, 0)=f_{4, 0}(x, 0)=0.
\end{align*}
From \eqref{A.4} and \eqref{A.5}, we derive the following additional boundary conditions for the viscous system \eqref{A.1}:
\begin{align}
\label{A.9}
f_3^{\varepsilon}(t, x, 0)|_{y=0}=f_4^{\varepsilon}(t, x, 0)|_{y=0}=0.
\end{align}
Similarly, from \eqref{A.6}, we have
\begin{equation}
f_3^{0}(t, x, 0)|_{y=0}=f_4^{0}(t, x, 0)|_{y=0}=0.\label{A.9-1}
\end{equation}
Therefore, the second row of the deformation matrix does not exhibit a strong boundary layer.

Now we follow Prandtl's approach to derive the boundary layer model for the viscoelastic equation \eqref{A.1}. We look for the solution to \eqref{A.1} in the following form:
\begin{equation}
\left\{
\begin{aligned}
\left(u^{\varepsilon}(t, x, y),v^{\varepsilon}(t, x, y)\right)&=\left(u^0(t, x, y),v^0(t, x, y)\right)+\left(u^b(t, x, \eta),\sqrt{\varepsilon}v^b(t, x, \eta)\right)+o(1),\\
\left(f_1^{\varepsilon}(t, x, y),f_3^{\varepsilon}(t, x, y)\right)&=\left(f_1^0(t, x, y),f_3^0(t, x, y)\right)+\left(f_1^b(t, x, \eta)+\sqrt{\varepsilon}f_3^b(t, x, \eta)\right)+o(1), \\
\left(f_2^{\varepsilon}(t, x, y),f_4^{\varepsilon}(t, x, y)\right)&=\left(f_2^0(t, x, y),f_4^0(t, x, y)\right)+\left(f_2^b(t, x, \eta), \sqrt{\varepsilon}f_4^b(t, x, \eta)\right)+o(1), \\
p^{\varepsilon}(t, x, y)&=p^0(t, x, y)+p^b(t, x, \eta)+o(1).
\label{A.16}
\end{aligned}\right.
\end{equation}
Here, $\eta=\frac{y}{\sqrt{\varepsilon}}$ represents the boundary layer variable and $(u^b, v^b, f_1^b, f_2^b, f_3^b, f_4^b, p^b)$ represents boundary layer corrections of the velocity field and deformation tensor which decay as $\eta\rightarrow +\infty$.

Substituting the ansatz \eqref{A.16} into \eqref{A.1} and collecting terms of the same order in $\varepsilon$, we obtain the leading-order outer flow $(u^0, v^0, f_1^0, f_2^0, f_3^0, f_4^0, p^0)$, which satisfies the ideal elastodynamic system \eqref{A.6}. The corresponding boundary layer profiles $(u,v,f_1,f_2,f_3,f_4,p)$ are then defined as follows:
\begin{align*}
\begin{cases}
u(t, x, \eta)=u^0(t, x, 0)+u^b(t, x, \eta),\\
v(t, x, \eta)=v^1(t, x, 0)+v^b(t, x, \eta)+\eta\partial_yv^0(t, x, 0),\\
f_1(t, x, \eta)=f_1^0(t, x, 0)+f_1^b(t, x, \eta),\\
f_2(t, x, \eta)=f_2^0(t, x, 0)+f_2^b(t, x, \eta),\\
f_3(t, x, \eta)=f_3^1(t, x, 0)+f_3^b(t, x, \eta)+\eta\partial_yf_3^0(t, x, 0),\\
f_4(t, x, \eta)=f_4^1(t, x, 0)+f_4^b(t, x, \eta)+\eta\partial_yf_4^0(t, x, 0),\\
p(t, x, \eta)=p^0(t, x, 0)+p^b(t, x, \eta),
\end{cases}
\end{align*}
where $(v^1,f_3^1,f_4^1)$ are corrections in the bulk at $O(\sqrt{\eps})$, which satisfies a linearized elastodynamic system.

By a straightforward computation, we derive the following equation for the boundary layer profile $(u, v, f_1, f_2, f_3, f_4, p)$:
\begin{align}
\label{A.19}
\begin{cases}
\partial_tu+u\partial_xu+v\partial_{\eta}u+\partial_xp-\partial_{\eta}^2u-f_1\partial_xf_1-f_3\partial_{\eta}f_1-f_2\partial_xf_2-f_4\partial_{\eta}f_2=0,\\
\partial_{\eta} p=0,\\
\partial_tf_1+u\partial_xf_1+v\partial_{\eta}f_1-f_1\partial_xu-f_3\partial_{\eta}u=0,\\
\partial_tf_2+u\partial_xf_2+v\partial_{\eta}f_2-f_2\partial_xu-f_4\partial_{\eta}u=0,\\
\partial_xu+\partial_{\eta}v=0,\quad \partial_xf_1+\partial_{\eta}f_3=0,\quad \partial_xf_2+\partial_{\eta}f_4=0.
\end{cases}
\end{align}

At the far field, the tangential velocity field $u$, the pressure $p$, and the first row of deformation matrix $(f_1,f_2)$ match the ideal flow in the following way:
\begin{equation}
\left\{
\begin{aligned}
\lim\limits_{\eta\to\infty}u(t, x, \eta)&=u^0(t, x, 0):=U(t, x),\\
p(t,x,\eta)&\equiv p^0(t,x,0):=P(t,x),\\
\lim\limits_{\eta\to\infty}f_1(t, x, \eta)&=f_1^0(t, x, 0):=F_1(t, x),\\
\lim\limits_{\eta\to\infty}f_2(t, x, \eta)&=f_2^0(t, x, 0):=F_2(t, x).\label{far}
\end{aligned}\right.
\end{equation}
Moreover, from \eqref{A.3} and \eqref{A.9}, we derive the following boundary conditions of $(u,v,f_3,f_4):$
\begin{align}
u(t, x, \eta)|_{\eta=0}=v(t, x, \eta)|_{\eta=0}=f_3(t, x, \eta)|_{\eta=0}=f_4(t, x, \eta)|_{\eta=0}=0.\label{A.21}
\end{align}
Combining \eqref{A.19}, \eqref{far} and \eqref{A.21} yields the boundary layer model \eqref{1.1}.

\bigbreak

\noindent{\bf Acknowledgements}.
The research of F. Xie is supported by the National Natural Science
Foundation of China under Grants 12671277 and 12271359. The research of Z. Zhang is supported by the General Research Fund from Hong Kong RGC (Project No.15300225) and Research Center for Nonlinear Analysis (Project No. P0046121) from The Hong Kong Polytechnic University.\\

\noindent{\bf Conflict of interest.} There is no conflict of interest.\\

\noindent{\bf AI-use disclosure.}  The authors declare that no AI model has been used in the paper.\\

\noindent{\bf Data availability.} This manuscript has no associated data.

\end{document}